\documentclass{amsart}

\usepackage{amsmath,amssymb,amsthm}
\usepackage{booktabs}
\usepackage{xcolor}
\usepackage[hyphens]{url}
\usepackage[colorlinks=true,linkcolor=blue!50!black,citecolor=blue!50!black,urlcolor=blue!50!black]{hyperref}

\DeclareFontFamily{U}{wncy}{}
\DeclareFontShape{U}{wncy}{m}{n}{<->wncyr10}{}
\DeclareSymbolFont{mcy}{U}{wncy}{m}{n}
\DeclareMathSymbol{\Sha}{\mathord}{mcy}{"58}

\newtheorem{theorem}{Theorem}[section]
\newtheorem{proposition}[theorem]{Proposition}
\newtheorem{lemma}[theorem]{Lemma}
\newtheorem{corollary}[theorem]{Corollary}
\newtheorem{conjecture}[theorem]{Conjecture}
\newtheorem{introtheorem}{Theorem}

\newtheorem*{introconjecture}{Conjecture}
\theoremstyle{definition}

\newtheorem{hypothesis}[theorem]{Hypotheses}
\theoremstyle{remark}
\newtheorem{remark}[theorem]{Remark}
\numberwithin{equation}{section}
\newcommand{\Q}{\mathbf{Q}}
\newcommand{\Z}{\mathbf{Z}}
\newcommand{\F}{\mathbf{F}}
\newcommand{\C}{\mathbf{C}}
\newcommand{\R}{\mathbf{R}}
\newcommand{\Zp}{\mathbf{Z}_p}
\newcommand{\Qp}{\mathbf{Q}_p}
\newcommand{\OK}{{\mathcal{O}_K}}
\newcommand{\fp}{\mathfrak{p}}
\newcommand{\fpb}{{\bar{\mathfrak{p}}}}

\newcommand{\ff}{\mathfrak{f}}
\newcommand{\iop}{{\iota_{\mathfrak{p}}}}
\newcommand{\psib}{{\bar{\psi}}}
\newcommand{\Shan}{{\Sha_{\mathrm{an}}}}

\newcommand{\Reg}{\operatorname{Reg}}
\newcommand{\ord}{\operatorname{ord}}
\newcommand{\adj}{\operatorname{adj}}
\newcommand{\tors}{\mathrm{tors}}
\newcommand{\Tam}{\textstyle\prod_v c_v}
\newcommand{\Frob}{\varpi}
\newcommand{\Lc}{\mathcal{L}}
\newcommand{\BDP}{\mathrm{BDP}}

\newcommand{\FourthDone}{301}

\newcommand{\FourthVthree}{122}
\newcommand{\FourthVtwo}{143}
\newcommand{\FourthVone}{30}
\newcommand{\FourthVzero}{6}
\newcommand{\FourthAnom}{90}
\newcommand{\FourthKs}{2, 4, 6}
\newcommand{\FourthShaFour}{291}
\newcommand{\FourthShaNine}{10}
\newcommand{\FourthNmin}{1.5\cdot10^{8}}
\newcommand{\FourthNmax}{2.9\cdot10^{12}}
\newcommand{\FourthAE}{\{4,8,9,16,18,32,9/2\}}

\newcommand{\FourthWhich}{\Shan\in\{4,9\}}

\begin{document}

\title[The $p$-adic BSD quotient in rank two]{Elliptic units, derived Heegner classes and the $p$-adic Birch and Swinnerton-Dyer quotient of CM elliptic curves of rank two}

\author{Miquel Calonge Fabregat}
\address{Universitat Rovira i Virgili, Tarragona, Spain}
\email{miquel1001@protonmail.com}

\subjclass[2020]{Primary 11G40, 11G05; Secondary 11G15, 11R23, 11Y40}
\keywords{Elliptic curves, complex multiplication, $p$-adic $L$-functions, $p$-adic heights, Heegner points, Birch and Swinnerton-Dyer conjecture}
\date{November 2026}

\begin{abstract}
Let $E\colon y^2=x^3-Dx$ have rank two and let $p\equiv1\pmod4$ be a prime of good reduction. We study the conjecture that the $p$-adic quotient $q_p$, the leading coefficient of the cyclotomic $p$-adic $L$-function divided by the Euler factor and the $p$-adic regulator, equals the complex quotient $q_\infty=L''(E,1)/(2\Omega_E\Reg_\infty)$, and that the latter is rational. In rank one the analogous statement is a theorem of Perrin-Riou and Rubin, because a Heegner point computes both sides; in rank two there is no such point. We show that $q_p$ is the constant of the quadratic term of the two-variable Katz measure; this is how the conjecture arises, and it computes $q_p$ modulo any power of $p$ from critical Hecke $L$-values. We then show that $q_p$ is the coordinate of two global classes: Rubin's rank-two elliptic-unit class, whose position in $E(\Q)\otimes\Qp$ is $q_p$ times an explicit vector, and the derived Heegner class over an auxiliary imaginary quadratic field $F$, for which we prove an exact identity. As consequences, whenever such a field $F$ satisfying our hypotheses exists, $q_p$ is a square in $\Qp$ times the product of the Tamagawa numbers of $E$, so that the unit $q_p/A_E$ of the $p$-adic Birch and Swinnerton-Dyer formula is a square when $\Sha$ is finite; if moreover the derived Heegner class has non-zero formal logarithm, the conjecture of Agboola and Castella for $E/F$ holds up to a $p$-adic unit, and with a theorem of Burungale, Castella and Skinner the anticyclotomic main conjecture holds integrally; and our conjecture is equivalent to the analytic form of theirs, with the products of anticyclotomic heights read sesquilinearly, which differs by a sign from its literal formulation. The conjecture is verified modulo $p^3$ for $5520$ pairs $(E,p)$ with conductors up to $3.9\cdot10^{12}$, including curves with analytic $\Sha$ of order $4$ and $9$.
\end{abstract}

\maketitle
\enlargethispage{4pt}

\section{Introduction}\label{sec:intro}

\subsection{The problem}
Let $E=E_D\colon y^2=x^3-Dx$, with $D$ a fourth-power-free integer, be an elliptic curve of rank two over $\Q$; it has complex multiplication by $\Z[i]$. Let $p\equiv1\pmod4$ be a prime of good reduction, so that $E$ is ordinary at $p$, and let $\alpha$ be the unit root of $X^2-a_pX+p$.
Two numbers measure the leading term of $L(E,s)$ at $s=1$. On the complex side,
\[q_\infty=\frac{L''(E,1)}{2\,\Omega_E\,\Reg_\infty}\in\R ,\]
where $\Omega_E$ is the real N\'eron period and $\Reg_\infty$ the N\'eron--Tate regulator of $E(\Q)$. On the $p$-adic side, if $c_2(p)$ is the coefficient of $T^2$ of the $p$-adic $L$-function of Mazur, Tate and Teitelbaum \cite{MTT} and $\Reg_\gamma$ the cyclotomic $p$-adic regulator (both normalised as in \S\ref{sec:katz}),
\[q_p=\frac{c_2(p)}{(1-\alpha^{-1})^2\,\Reg_\gamma}\in\Qp .\]
The Birch and Swinnerton-Dyer conjecture and its $p$-adic analogue \cite{MTT} predict that both are equal to $\#\Sha(E/\Q)\Tam/\#E(\Q)_\tors^2$. We study the statement that does not involve $\Sha$:

\begin{introconjecture}
$q_\infty$ is a rational number, and $q_p=q_\infty$ in $\Qp$.
\end{introconjecture}

In rank zero this follows from the interpolation property of the $p$-adic $L$-function. In rank one the corresponding identity is a theorem of Perrin-Riou \cite{PR87} and Rubin \cite[Thm.~8.2]{Rubin92}, and the reason is a single global object. A Heegner point computes the complex derivative, by the formula of Gross and Zagier \cite{GZ86}, and the $p$-adic derivative, by its $p$-adic analogue; Rubin obtains the $p$-adic side for CM curves from elliptic units. Both quotients are then the same ratio of heights of the same point, and no arithmetic invariant enters.
In rank two there is no such point. Over $K=\Q(i)$ the Heegner hypothesis fails at $2$, and over an auxiliary field $F$ that satisfies it the Heegner point is torsion, because $L(E/F,s)$ vanishes to order three. Moreover, $\Sha$ is not known to be finite for a single curve of analytic rank at least two \cite{CLS1}, and the Iwasawa main conjecture, an equality of ideals, determines $q_p$ only up to a $p$-adic unit.

The thesis of this paper is that in rank two $q_p$ is nevertheless attached to global objects, in two ways. First, $q_p$ is the constant of the quadratic term of the two-variable Katz measure; this is how the conjecture was found, and it makes $q_p$ computable to any $p$-adic precision from critical Hecke $L$-values. Second, $q_p$ is the coordinate of two canonical global classes: Rubin's rank-two elliptic-unit class over $K=\Q(i)$, and the derived Heegner class over an auxiliary imaginary quadratic field $F$. From the second class we extract what can be proved: an exact identity, the square class of $q_p$, and the equivalence of the conjecture with the conjecture of Agboola and Castella for $E/F$.

\subsection{Results}
Let $H$ be the matrix of cyclotomic $p$-adic heights of generators $P_1,P_2$ of $E(\Q)$ modulo torsion, $j$ the anticyclotomic height pairing of $[i]P_1$ and $P_2$ (\S\ref{ss:heights}), and $\det H'=p^{-2}\det H$, $j'=p^{-1}j\bmod p$.
Let $M_k$ be the central value of $L(\psib^{\,2k(p-1)+1},s)$, normalised as in \S\ref{sec:setup}, and $M=M_1$.

\emph{The Katz measure.} Perrin-Riou's leading-term formula \cite{PR84} describes the quadratic term of the two-variable $p$-adic $L$-function of $E$ over $K$ up to a unit. We identify the constant.

\begin{introtheorem}[Theorem~\ref{thm:three}, Corollary~\ref{prop:Cp}, Theorem~\ref{thm:digits}]\label{thm:A}
Assume $\Sha(E/\Q(i))[p^\infty]$ finite and $\det H\ne0$. The quadratic term of the Katz measure at $\psi_E$ is $q_p$ times the explicit quadratic form
\[\tfrac14\,\iop(f_0\omega_E)\,\Omega_p\,e^2(1-\alpha^{-1})^2\bigl[(a-b)^2\det H-(a+b)^2j^2\bigr]\qquad(e=p-1).\]
If moreover $a_p\not\equiv1\pmod p$, then with $c=\iop(f_0\omega_E)(a_p-1)^2\in\F_p^\times$
\[M'\equiv-c\,q_p\,j'^2\pmod p,\]
and similarly in the cyclotomic and $\fp$-directions; and if $v_p(j)=1$, the values $M_1,\dots,M_m$ determine $q_p$ modulo $p^{m}$ $(p\ge m+3)$ by an explicit formula.
\end{introtheorem}

The three directions of Perrin-Riou's pencil, cyclotomic, anticyclotomic and $\fp$-adic, carry the same constant, and it is $q_p$. Comparing the values of $q_p$ so obtained with $q_\infty$ led to the conjecture, and the second part of the theorem is the algorithm behind the evidence below.

\emph{Elliptic units.} Rubin \cite{Rubin92} constructs from elliptic units a rank-two class $x^*$ in $E(K)\otimes\Qp$ and proves exact formulas for its pairings and logarithms in terms of the Katz measure. Combined with Theorem~\ref{thm:A}, they locate the class.

\begin{introtheorem}[Theorems~\ref{thm:position}, \ref{thm:fourth}]\label{thm:B}
Assume $\Sha(E/\Q(i))[p^\infty]$ finite, $\det H\ne0$ and $\Delta=\det H-j^2\ne0$; $p$ may be anomalous. Then Rubin's class $x^*$ is, in an explicit basis of $E(\Q)\otimes\Qp$, equal to $q_p$ times the vector $\tfrac{1}{2p}\iop(f_0\omega_E)(1-\alpha/p)\adj(H+\iop(i)^{-1}J)\ell$, where $\ell$ is the vector of formal logarithms of $P_1,P_2$ and $J=\left(\begin{smallmatrix}0&j\\-j&0\end{smallmatrix}\right)$. Consequently the linear term $\kappa_1$ of the Katz measure on the branch through the conjugate character satisfies
\[\frac{\kappa_1}{p-1}=\tfrac12\,\iop\Bigl(\frac{\bar\omega_E}{f_0}\Bigr)\Bigl(1-\frac\alpha p\Bigr)^2q_p\,h^+(Q_\lambda),\]
where $Q_\lambda$ spans the line of vanishing logarithm.
\end{introtheorem}

Thus the conjecture says that the position of Rubin's class in $E(\Q)\otimes\Qp$ is the one given by the rational number $q_\infty$. This is the rank-two analogue of Rubin's rank-one formula, which computes a rational point from the $p$-adic $L$-function; what is missing in rank two is an archimedean formula for the class.

\emph{Derived Heegner classes.} Now let $F$ be an imaginary quadratic field of odd discriminant in which $p$ and every prime dividing $N_E$ split, such that the twist $E^\chi$ of $E$ by $F$ has analytic rank one; then $E(F)$ has rank three, Heegner points over $F$ exist, and those defined over $F$ itself are torsion. Howard's $\Lambda$-adic Heegner class over the anticyclotomic $\Zp$-extension of $F$ therefore vanishes at the origin, and we let $\mathcal D$ be its derivative there. Howard's $\Lambda$-adic Gross--Zagier theorem \cite{Howard05}, Hsieh's main conjecture for CM fields \cite{Hsieh14} and Nekov\'a\v{r}'s Selmer complexes \cite{Nekovar06} relate $\mathcal D$ to the same quadratic structure, and the result is an identity of elements, not of ideals.

\begin{introtheorem}[Theorems~\ref{thm:heegner}, \ref{thm:squares}, \ref{thm:tamsquare}]\label{thm:C}
Under the hypotheses of \S\ref{ss:hyp}, $\mathcal D=\zeta\,Q_v$ with $\zeta\in\Qp^\times$, where $Q_v=h_a(P_2,P_3)P_1-h_a(P_1,P_3)P_2$ is the vector of $E(\Q)\otimes\Qp$ orthogonal to a generator $P_3$ of $E^\chi(\Q)$ for the anticyclotomic height $h_a$ of $F$, and
\[q_p=\frac{1}{S_F}\Bigl(\frac{2\,\ell_a(\sigma_a)\,\zeta}{(1-\alpha^{-1})^{2}}\Bigr)^{2},\]
with $S_F\in\Q^\times$ explicit in terms of the rank-one quotient of $E^\chi$; after normalising by the Manin constant, $s_F=c_E^{-2}S_F$ equals $q_\infty(E^\chi)$ or $\tfrac12q_\infty(E^\chi)$ (Lemma~\ref{lem:kappa}). Consequently $q_p\,s_F\in(\Qp^\times)^2$, and in fact $q_p\in\Tam(E)\cdot(\Qp^\times)^2$.
\end{introtheorem}

The last statement combines the identity with a theorem on the rational number $s_F$: comparing $E^\chi$ with a real twist of rank zero through a genus character, the explicit Gross--Zagier formula of Cai, Shu and Tian \cite{CST14} and the Birch and Swinnerton-Dyer formula in analytic rank zero for CM curves of Burungale and Flach \cite{BF24} show that $\Tam(E)\,s_F$ is a rational square (Theorem~\ref{thm:tamsquare}). The square class is unconditional given the hypotheses, which require one auxiliary field $F$ for $(E,p)$ (Remark~\ref{rem:existsF}), and it is the first information on the unit $q_p/A_E$ beyond its valuation: if $\Sha(E/\Q)$ is finite, $q_p/A_E$ is a square in $\Zp^\times$. The identity carries no undetermined sign: the two sign conventions in the comparison of Howard's height with ours cancel (Lemma~\ref{lem:howard-sign}). The archimedean counterpart of the identity is the Birch and Swinnerton-Dyer formula over $F$, which expresses $q_\infty s_F$ through the analytic order of $\Sha(E/F)$ (Theorem~\ref{thm:bsdF}).

\emph{The conjecture of Agboola and Castella.} The derived Heegner class gives the leading term of the Bertolini--Darmon--Prasanna $p$-adic $L$-function $L^\BDP_p$ of $E/F$, by the explicit reciprocity law of Castella and Hsieh \cite[Thm.~5.7]{CastellaHsieh18} together with Castella's big logarithm \cite[Thm.~3.7]{Castella20}, and Agboola and Castella \cite{AC21} have conjectured a formula for that leading coefficient in terms of a regulator of anticyclotomic heights.

\begin{introtheorem}[Theorems~\ref{thm:AC}, \ref{thm:bdpMC}, Corollary~\ref{cor:ACunit}]\label{thm:D}
Assume moreover that $Q_v$ has non-zero formal logarithm. Then $L^\BDP_p$ vanishes to order exactly two, and:
\begin{enumerate}
\item the conjecture of Agboola and Castella \cite[Conj.~1.1, 4.2]{AC21} holds for $(E,F,p)$ up to a $p$-adic unit;
\item its exact form, with $\#\Sha(E/F)$ replaced by its analytic value and the products of anticyclotomic heights read sesquilinearly, holds if and only if $q_\infty\in\Q$ and $q_p=q_\infty$, that is, if and only if the conjecture above holds;
\item the square class of the leading coefficient is the one predicted by the conjecture when $\Shan(E/F)$ is a square in $\Qp$, as it is when $\Sha(E/F)$ is finite and equal to $\Shan(E/F)$;
\item the anticyclotomic Iwasawa--Greenberg main conjecture for $(E,F,p)$, $\mathrm{char}(X_{\rm Gr})=(L^\BDP_p)$, holds integrally.
\end{enumerate}
\end{introtheorem}

Part (1) needs no main conjecture for $L^\BDP_p$: it follows from the exact identity and from the known $p$-parts of the Birch and Swinnerton-Dyer formula for $E$ and $E^\chi$. Part (4) combines (1) with the main conjecture after inverting $p$, which Burungale, Castella and Skinner prove by base change without excluding complex multiplication \cite[Thm.~1.2.4(a)]{BCS25}; for these residually dihedral primes we check the irreducibility conditions of their proof. The sesquilinear reading in (2) is the one dictated by the $\Lambda$-adic theory, in which heights are sesquilinear. Read literally, with the products in the ring, the conjecture for $L^\BDP_p=\Lc_p^2$ differs by a sign from the square of the conjecture of Agboola and Castella for $\Lc_p$ \cite[Lemma~3.14, Conj.~4.3]{AC21}, and it is numerically incompatible with the two examples we computed; the formulation of Castella, Hsu, Kundu, Lee and Liu \cite{CHKLL} is compatible with the sesquilinear sign for a natural choice of bases, although its results, which hold up to units, do not determine it (\S\S\ref{sec:AC}, \ref{ss:readings}).

\subsection{Evidence}
The quotient $q_p$ is computable from critical values of Hecke $L$-functions and from $p$-adic heights, with an algorithm whose cost grows like $\sqrt{N_E}$, so that conductors far beyond the reach of modular symbols can be treated. For $5941$ pairs $(E_D,p)$ with $p\le41$ and conductors up to $3.9\cdot10^{12}$ we find $q_p\equiv q_\infty\pmod{p^3}$ in all the $5520$ pairs where Theorem~\ref{thm:A} applies; $q_\infty$ involves $\Shan=4$ for $226$ curves and $\Shan=9$ for $20$, recognised numerically as integers; these are prime to $p$ and invisible to the main conjecture, but the $p$-adic values detect them modulo $p^3$. For $E_{56}$ one has $q_5\equiv1\pmod{5^{12}}$. The identity of Theorem~\ref{thm:B} with $q_\infty$ holds to five $5$-adic digits for three curves and shows no discrepancy in any computed digit for $301$ curves with $\Shan\in\{4,9\}$ ($295$ of them with at least one digit). The identity of Theorem~\ref{thm:C} is verified from Heegner points of conductor $25$ and $125$ for $E_{56}$ over $\Q(\sqrt{-31})$, modulo $25$, and blindly modulo $5$ for $E_{-289}$ over $\Q(\sqrt{-111})$, sign included, and the archimedean identity for $79$ pairs $(E,F)$ (\S\ref{sec:evidence}).

\subsection{Relation to previous work}
The quadratic structure behind Theorem~\ref{thm:A} is Perrin-Riou's \cite{PR84}; the cyclotomic second derivative was expressed through critical values by Banwait \cite{Banwait}, whose bracket is one of the three directions. Theorem~\ref{thm:B} combines Rubin's Theorem~2.5, Proposition~9.4 and Theorem~9.5 \cite{Rubin92} with the exact constant of Theorem~\ref{thm:A}; Castella \cite{Castella26} shows that the generalised Kato class spans the line of vanishing logarithm, and together with Rubin's class the Selmer group. The conjecture is the analytic form of the Generalized Perrin-Riou Conjecture of Burns, Kurihara and Sano \cite{BKS} for these curves (Remark~\ref{rem:BKS}). Theorem~\ref{thm:C} rests on Howard's $\Lambda$-adic Gross--Zagier theorem \cite{Howard05}, Hsieh's main conjecture \cite{Hsieh14}, Nekov\'a\v{r}'s Selmer complexes \cite{Nekovar06} and Rubin's rank-one theorem \cite{Rubin92}; Theorem~\ref{thm:D} on the explicit reciprocity law of Castella and Hsieh \cite[Thm.~5.7]{CastellaHsieh18}, on Castella's big logarithm \cite[Thm.~3.7]{Castella20} and on \cite{AC21,BDP13}. Agboola and Castella prove their formula up to a $p$-adic unit for the characteristic power series under a large-image hypothesis that excludes complex multiplication \cite[Thm.~1.3]{AC21}; Castella, Hsu, Kundu, Lee and Liu prove the algebraic formula up to a unit under hypotheses that allow it \cite[Conj.~1.1, Thm.~5.1]{CHKLL}, and Sano deduces the conjecture up to a unit from the integral anticyclotomic main conjecture \cite[Thm.~4.13]{Sano23}. In our setting Burungale, Castella and Skinner prove the main conjecture after inverting $p$ \cite[Thm.~1.2.4(a)]{BCS25}, which with Sano's theorem gives the formula only up to a power of $p$. Theorem~\ref{thm:D} gives the formula up to a unit directly, and with it the integral main conjecture (Theorem~\ref{thm:bdpMC}); it also expresses the unit through $q_p$: when $q_\infty\in\Q$ and $\Sha(E/F)$ is finite, it is $q_p\,\Shan(E/F)/(q_\infty\,\#\Sha(E/F))$. The numerical evidence for the cyclotomic $p$-adic Birch and Swinnerton-Dyer conjecture in rank two \cite{MTT,SteinWuthrich}, and for its analogue for abelian surfaces \cite{BMS}, is limited to small conductors; Bernardi, Goldstein and Stephens \cite{BGS} tested the $\fp$-adic version for CM curves, including $E_{17}$ at $p=5$.

\subsection{Organisation}
Section~\ref{sec:setup} fixes the setting. Section~\ref{sec:katz} proves Theorem~\ref{thm:A} and states the conjecture; Section~\ref{sec:rubin} proves Theorem~\ref{thm:B}; Sections~\ref{sec:heegner} and~\ref{sec:squares} prove Theorem~\ref{thm:C} and the Birch and Swinnerton-Dyer formula over $F$; Section~\ref{sec:AC} proves Theorem~\ref{thm:D}. All the numerical evidence is in Section~\ref{sec:evidence}, and Section~\ref{sec:final} discusses what remains open.

\subsection{Reproducibility}
All computations were done with PARI/GP \cite{PARI}, with independent checks in SageMath \cite{Sage}. The large-scale computations were run on rented cloud machines (vast.ai), the others on a personal computer. Code, job scripts and the outputs of the computations are provided as supplementary material, archived at Zenodo, \url{https://doi.org/10.5281/zenodo.23982770}. Where an original log was not preserved (the rank-one calibration of \S\ref{ss:rankone}, and the twelve-value and modular-symbol checks for $E_{56}$), the computation was rerun in November 2026 and the new output is archived instead.

\section{Setting}\label{sec:setup}

\subsection{Curves and primes}
Let $D\ne0$ be fourth-power-free and $E=E_D\colon y^2=x^3-Dx$. It has complex multiplication by $\OK=\Z[i]$, $K=\Q(i)$, with $[i](x,y)=(-x,iy)$; the model is global minimal (an elementary consequence of Kraus's criterion \cite[Prop.~2 and \S D.1]{Kraus89}). Let $\psi=\psi_E$ be the Gr\"ossencharacter of $E$, of infinity type $(1,0)$ and conductor $\ff=(f_0)$, so that $L(E,s)=L(\psi,s)$, $N_E=4\,\mathrm N\ff$ and
\begin{equation}\label{eq:f0}
f_0=(1+i)^{e_f}\sqrt{N_E^{\mathrm{odd}}},\qquad e_f=v_2(N_E)-2 ;
\end{equation}
at an odd prime $q\mid D$ the reduction is additive and the $q$-part of $\ff$ is $(q)$. The lattice of $\omega=dx/2y$ is $\Omega_\infty\OK$ with $\Omega_\infty=\omega_1$ if $D>0$ and $\Omega_\infty=\tfrac{1+i}{2}\omega_1$ if $D<0$, $\omega_1$ the least positive real period, and the real N\'eron period is $\Omega_E=\omega_E\Omega_\infty$ with $\omega_E=2$ if $D>0$ and $\omega_E=1-i$ if $D<0$. Replacing $f_0$ or $\Omega_\infty$ by a unit multiple $u$ changes the quantities below coherently: for each weight that we use, $\mathcal N(k,l)$, the measure and the factors such as $\iop(f_0\omega_E)$ change by compatible powers of $u$, so that both sides of each identity change by the same fourth root of unity, and $q_p$ does not change.
We assume throughout that $E(\Q)$ has rank two and root number $+1$, and fix points $P_1,P_2$ generating $E(\Q)$ modulo torsion.

Let $p\equiv1\pmod4$ be a prime with $p\nmid N_E$. Write $a_p=2a_0$, $p=a_0^2+b_0^2$ with $b_0>0$, and $\fp=(\Frob)$ with $\Frob=\psi(\fp)=a_0+b_0i$. Let $\iop\colon K\to\Qp$ be the embedding with $\iop(\Frob)\in p\Zp$, so that $\alpha=\iop(\bar\Frob)$ is the unit root of $X^2-a_pX+p$ and $\beta=p/\alpha$. We write $W$ for the ring of integers of the completion of the maximal unramified extension of $\Qp$, and $e=p-1$. The prime $p$ is \emph{anomalous} if $a_p\equiv1\pmod p$, that is, if $1-\alpha^{-1}$ is not a unit.

\subsection{The Katz measure and critical values}
For integers $k,l\ge0$ put $w=k+l+1$ and
\begin{equation}\label{eq:EK}
\mathcal N(k,l)=\Bigl(\frac{f_0}{\Omega_\infty}\Bigr)^{w}\Bigl(1-\frac{\Frob^{w}}{p^{l+1}}\Bigr)\Bigl(1-\frac{\Frob^{w}}{p^{k+1}}\Bigr)\,k!\,\pi^{l}\,L(\psib^{\,w},k+1),
\end{equation}
which lies in $K$ by Damerell's theorem and is read in $\Qp$ through $\iop$. Let $\mu$ be the measure of Bannai and Kobayashi \cite[\S3]{BK10} on $G_\infty=\mathrm{Gal}(K(\ff p^\infty)/K)$, a normalisation of the measure of Katz \cite{Katz76}, attached to the lattice $\Omega_\infty\OK$; let $x,y\colon G_\infty\to\Zp^\times$ be the $p$-adic avatars of $\psi$ and $\psib$ (the characters giving the action on $E[\fp^\infty]$ and $E[\fpb^\infty]$), and $\mu_\psi=x\cdot\mu$.

\begin{proposition}[Bannai--Kobayashi]\label{prop:moments}
For integers $k\ge l\ge0$ with $w=k+l+1$ odd,
\begin{equation}\label{eq:moments}
\int_{G_\infty}x^{k}y^{-l}\,d\mu_\psi=(-1)^{k+l}\,\Omega_p^{\,w}\,\iop\bigl(\mathcal N(k,l)\bigr),
\end{equation}
where $\Omega_p=\theta'(0)^{-1}\in W^\times$ for an isomorphism of formal groups $\theta\colon\widehat E\to\widehat{\mathbf G}_m$ over $W$, i.e.\ $\theta(t)=\exp(\Omega_p^{-1}\log_{\widehat E}(t))-1$ with $\log_{\widehat E}$ the formal logarithm \cite[\S3.1]{BK10}. For odd $w$ the finite part of $\psib^{\,w}$ is that of $\psib^{\pm1}$, whose values are fourth roots of unity, so $\psib^{\,w}$ has conductor $\ff$ and $L(\psib^{\,w},s)$ is the primitive $L$-function; all the weights used below are odd.
\end{proposition}
\begin{proof}
This is \cite[Thm.~3.11]{BK10} for the Hecke character $\psi^{k+1}\psib^{-l}$, of infinity type $(k+1,-l)$ with trivial local character at $\fp$: the Euler factor there is $1-\Frob^w/p^{l+1}$, removing the factor at $\fpb$ gives the second one, $2\pi/\sqrt{d_K}=\pi$, and $L(\psi^{-k-1}\psib^{\,l},0)=L(\psib^{\,w},k+1)$ because $\psi^{-1}=\psib\,\mathrm N^{-1}$.
\end{proof}

\begin{lemma}\label{lem:period}
If $\varphi_W$ is the Frobenius of $W$, then $\Omega_p^{\varphi_W}=\alpha^{-1}\Omega_p$; in particular $\Omega_p^{\,p-1}\equiv a_p^{-1}\pmod{\mathfrak m_W}$.
\end{lemma}
\begin{proof}
$\widehat E$ is a Lubin--Tate formal group over $\Zp$ for the uniformiser $\iop(\Frob)$ and $\widehat{\mathbf G}_m$ for $p=\alpha\,\iop(\Frob)$, so by \cite[Thm.~1, Lemma~2]{LT65} there is an isomorphism $\theta$ over $W$ with $\theta^{\varphi_W}=\theta\circ[\alpha]$; any other differs by an automorphism of $\widehat{\mathbf G}_m$ over $\Zp$, which does not change $\theta'(0)^{\varphi_W}/\theta'(0)=\alpha$. See also \cite[II.4.3]{deShalit}.
\end{proof}

The values that enter below are $T_1=\mathcal N(2e,0)$, $M=\mathcal N(e,e)$, $T_3=\mathcal N(0,2e)$, all for $L(\psib^{2p-1},s)$, and more generally the central values
\[M_k=\mathcal N(ke,ke)\qquad(k\ge1),\]
of $L(\psib^{\,2ke+1},s)$ at $s=ke+1$, whose root number is $+1$. The combination $B=T_1-2M+T_3$ is the bracket of Banwait \cite[Cor.~6.8]{Banwait}, and, for $p$ non-anomalous, $T_1$ is congruent modulo $\fp^3$ to a unit multiple of the quantity $c_p^+(E)$ of Coates, Liang and Sujatha \cite[Thm.~2.2]{CLS1}.

\subsection{Heights}\label{ss:heights}
Let $\sigma(t)$ be the $p$-adic sigma function of Mazur and Tate \cite[Thm.~3.1]{MT91} in $t=-x/y$, computed as in \cite[\S3, Algorithm~3.1]{MST}. For $K=\Q(i)$ one has $s_2=0$ \cite[p.~54]{PR84}, so $\sigma$ is the formal Weierstrass sigma function; it lies in $\Q[[t]]$ and does not depend on $p$. For $S\in E(K)$ in the formal group at $\fp$ and $\fpb$ and in the identity component at the bad places, write $x(S)=A/\delta_S^2$ with $A,\delta_S\in\OK$ coprime and put
\begin{equation}\label{eq:hp}
h_\fp(S)=-\log_p\iop\bigl(\sigma(t_S)/\delta_S\bigr),\qquad h_\fpb(S)=h_\fp(\bar S),
\end{equation}
extended to quadratic forms on $E(K)\otimes\Q$; these are the heights of \cite[Ch.~III]{PR84} attached to the two projections of $\OK\otimes\Zp=\Zp\times\Zp$. The cyclotomic and anticyclotomic heights over $K$ are $h^\pm=h_\fp\pm h_\fpb$, with bilinear forms $\langle\ ,\ \rangle^\pm$. We put
\[H=\bigl(\langle P_a,P_b\rangle^{+}\bigr)_{a,b=1,2},\qquad j=\langle[i]P_1,P_2\rangle^{-}.\]
On $E(\Q)$ one has $h^+(P)=-\log_p\bigl(\sigma(t_P)^2/d(P)\bigr)$, $d(P)$ the denominator of $x(P)$; this is the negative of the pairing $\langle P,P\rangle_\lambda=\log_p(\sigma^2/d)$ of \cite[\S II.5, p.~33]{MTT}. In rank two this does not change $\det H$, so $\det H$ is the cyclotomic $p$-adic regulator of \cite{MTT}; in rank one the sign matters, and we use $h^+$ (see the discussion after \eqref{eq:c1}). It is also $-2p$ times the height $h_p$ of \cite[(1.1)]{MST}, and it is the value returned by \texttt{ellpadicregulator} in PARI/GP, which we checked to $O(5^{15})$ for $E_{56}$. For $p$ non-anomalous one has $v_p(\det H)\ge2$ and $v_p(j)\ge1$; we put $\det H'=p^{-2}\det H$ and $j'=p^{-1}j$ modulo $p$. Since $\overline{[i]P_1}=-[i]P_1$, $j=2\langle[i]P_1,P_2\rangle_\fp$, and $j\ne0$ by a theorem of Bertrand \cite[Cor.~3]{Bertrand84}, as $[i]P_1$ and $P_2$ are $K$-linearly independent.

\section{\texorpdfstring{The quotient $q_p$ and the Katz measure}{The quotient q\_p and the Katz measure}}\label{sec:katz}\label{sec:constant}\label{sec:quadratic}

\subsection{\texorpdfstring{The $p$-adic quotient}{The p-adic quotient}}
Let $L_p(E,T)$ be the cyclotomic $p$-adic $L$-function of Mazur, Tate and Teitelbaum \cite{MTT} divided by the real period $\Omega_E=\Omega_E^+$, in the variable $T=\gamma^{s-1}-1$, $\gamma=1+p$, and let $c_2(p)$ be its coefficient of $T^2$. Since MTT differentiate in $s$ and use $\log_p$ itself in the regulator, their conjecture \cite[p.~38]{MTT} reads $c_2(p)\log_p(\gamma)^2=(1-\alpha^{-1})^2A_E\det H$ in rank two. Put
\begin{equation}\label{eq:qp}
\Reg_\gamma=\frac{\det H}{\log_p(1+p)^2},\quad q_p=\frac{c_2(p)}{(1-\alpha^{-1})^2\,\Reg_\gamma},\quad A_E=\frac{\#\Sha(E/\Q)\Tam}{\#E(\Q)_\tors^2}.
\end{equation}
The $p$-adic Birch and Swinnerton-Dyer conjecture of \cite{MTT} predicts $q_p=A_E$. By the cyclotomic main conjecture \cite[Thm.~12.3 and p.~67]{Rubin91}, $L_p(E,T)$ generates the characteristic ideal of the $p^\infty$-Selmer group over the cyclotomic $\Zp$-extension. By Perrin-Riou's theorem \cite[Thm.~V.8]{PR84}, applied to the cyclotomic character, its order of vanishing is at least $2+r$, where $r$ is the $\Zp$-corank of $\Sha(E/K)[\fp^\infty]$, and it equals $2$ if and only if $r=0$ and the cyclotomic height is non-degenerate, i.e.\ $\det H\ne0$ (Lemma~\ref{lem:dictionary} below). Hence $c_2(p)\ne0$ if and only if $\Sha(E/\Q)[p^\infty]$ is finite and $\det H\ne0$, and then $v_p(q_p)=v_p(\#\Sha(E/\Q)[p^\infty])$ by Schneider's leading-term formula \cite{Schneider85} or \cite[Thm.~V.8(3)]{PR84}, the Tamagawa numbers and the torsion being prime to $p$; if moreover $\Sha(E/\Q)$ is finite, so that $A_E$ is defined, $q_p/A_E\in\Zp^\times$; to know that $c_2(p)\ne0$ one needs $c_2(p)$ modulo $p^{v_p(c_2(p))+1}$, which modular symbols provide for small conductors.

\begin{lemma}[{\cite[Cor.~3.12]{BK10}}]\label{lem:BK}
Let $L^{\mathrm{Katz}}_p\in W[[T]]$ be defined by $L^{\mathrm{Katz}}_p((1+p)^{s}-1)=\int\langle xy\rangle^{s}\,d\mu_\psi$, where $\langle\ \rangle$ is the projection to $1+p\Zp$. Then $L_p(E,T)=(f_0\,\omega_E\,\Omega_p)^{-1}(1+T)^{m}L^{\mathrm{Katz}}_p(T)$ for some $m\in\Zp$.
\end{lemma}
\begin{proof}
For a Dirichlet character $\chi$ of $p$-power conductor, \cite[Cor.~3.12]{BK10} gives $\Omega_p^{-1}\int x(\chi\circ\mathrm N)\,d\mu=\nu\cdot\mathcal L^{MS}(\chi)$ with $\mathcal L^{MS}$ the Mazur--Swinnerton-Dyer measure and $\Omega_E=\nu\,\Omega_\infty/f_0$, i.e.\ $\nu=f_0\omega_E$; this is the interpolation property of $L_p(E,T)$ \cite[\S I.14]{MTT}, and a change of the system of roots of unity multiplies the measure by a group element, whose power series is $(1+T)^m$.
\end{proof}

\subsection{The quadratic term of the Katz measure}
The characters $x^{e}$ and $y^{e}$ are trivial on the torsion of $G_\infty$, so the moments $\int x^{ae}y^{-be}\,d\mu_\psi$ are the values at $(S,T)=(u_1^{a}-1,u_2^{-b}-1)$ of a power series $G(S,T)\in W[[S,T]]$, where $u_1,u_2\in1+p\Zp$ are the values of $x^e,y^e$ on topological generators. Put
\[\mathcal F(a,b)=G(u_1^{a}-1,\,u_2^{-b}-1)\qquad(a,b\in\Zp),\]
and let $\mathcal F=\sum_{n\ge0}\mathcal F_n$ be its expansion in homogeneous polynomials. Since $u_i^{a}-1=\sum_{l\ge1}(\log_pu_i)^la^l/l!$ with $v_p(\log_pu_i)=1$, the coefficients of $\mathcal F_n$ lie in $p^{\,n-v_p(n!)}W$.

\begin{proposition}\label{prop:quadratic}
Assume $\Sha(E/K)[\fp^\infty]$ finite. Then $\mathcal F_0=\mathcal F_1=0$; consequently $\int x^{ae}y^{-be}\,d\mu_\psi\equiv \mathcal F_2(a,b)\pmod{p^3W}$ for $a,b\in\Z_{\ge0}$.
\end{proposition}
\begin{proof}
By the two-variable main conjecture, as in the proof of Theorem~\ref{thm:three}(1) below, the Katz measure is a unit of $W[[S,T]]$ times a characteristic power series $\mathcal H_\fp$ of the dual $\fp^\infty$-Selmer group over $K_\infty$, and by Perrin-Riou's Theorem~\ref{thm:PR}, with $n=\operatorname{rank}_{\OK}E(K)=2$, $\mathcal H_\fp$ vanishes to order at least two at the origin. This kills the terms of degree $\le1$; terms of degree $n\ge3$ lie in $p^3W$.
\end{proof}

\subsection{Perrin-Riou's pencil}
Let $K_\infty=K(E[p^\infty])$ and $\kappa,\kappa^*$ the characters of $\mathrm{Gal}(K_\infty/K)$ giving the action on $E[\fp^\infty]$ and $E[\fpb^\infty]$. To $\rho=\kappa^{s}\kappa^{*s^*}$ Perrin-Riou attaches a $p$-adic height pairing between the two eigencomponents of $E(K)\otimes\Zp$, and proves:

\begin{theorem}[Perrin-Riou {\cite[Thm.~V.17]{PR84}}]\label{thm:PR}
Let $\mathcal H_\fp$ be a characteristic power series of the dual of the $\fp^\infty$-Selmer group of $E$ over $K_\infty$, $n=\operatorname{rank}_{\OK}E(K)$, and assume $\Sha(E/K)[\fp^\infty]$ finite. Then $\mathcal H_\fp(\kappa^{s}\kappa^{*s^*})$ is congruent modulo $(s,s^*)^{n+1}$ to a $p$-adic unit times explicit Euler factors, $\#\Sha(E/K)[\fp^\infty]$ and the discriminant of $s\langle\ ,\ \rangle_{\kappa,\fp}+s^*\langle\ ,\ \rangle_{\kappa^*,\fp}$.
\end{theorem}

For $E$ defined over $\Q$ the discriminant is computed in \cite[pp.~118--119]{PR84}: if $t,t^*$ are the two projections and $M_0=(\langle t(P_a),t^*(P_b)\rangle_{\kappa,\fp})_{a,b}$, then complex conjugation shows that the second pairing has matrix ${}^tM_0$, and the discriminant is $\det(sM_0+s^*\,{}^tM_0)$.

\begin{lemma}\label{lem:dictionary}
Let $r=\iop(i)$ and $J=\left(\begin{smallmatrix}0&j\\-j&0\end{smallmatrix}\right)$. Then $M_0=\tfrac14(H+r^{-1}J)$ and
\begin{equation}\label{eq:Dab}
\det(sM_0+s^*\,{}^tM_0)=\tfrac1{16}\bigl[(s+s^*)^2\det H-(s-s^*)^2j^2\bigr].
\end{equation}
\end{lemma}
\begin{proof}
The idempotent of $\OK\otimes\Zp$ mapping to $(1,0)$ under $(\iop,\iota_\fpb)$ is $\tfrac12(1\otimes1+i\otimes r^{-1})$, so $t(x)=\tfrac12(x+r^{-1}[i]x)$ and $t^*(x)=\tfrac12(x-r^{-1}[i]x)$. For $b=\langle\ ,\ \rangle_\fp$ one has $b(\gamma x,y)=b(x,\bar\gamma y)$, hence $b([i]x,[i]y)=b(x,y)$, $b(x,[i]y)=-b([i]x,y)$ and, as $r^{-2}=-1$, $b(t(x),t^*(y))=\tfrac12[b(x,y)+r^{-1}b([i]x,y)]$. For rational points $h^+=2h_\fp$, so $b(P_a,P_b)=\tfrac12H_{ab}$, and $(b([i]P_a,P_b))$ is antisymmetric with entry $\tfrac12j$. This gives $M_0$; then $sM_0+s^*\,{}^tM_0=\tfrac14((s+s^*)H+(s-s^*)r^{-1}J)$, and $\det(aH+bJ)=a^2\det H+b^2j^2$.
\end{proof}

The cyclotomic character is $\kappa\kappa^*$, i.e.\ $(s,s^*)=(1,1)$, where \eqref{eq:Dab} is $\tfrac14\det H$; the anticyclotomic direction $(1,-1)$ gives $-\tfrac14j^2$, in agreement with Perrin-Riou's remark \cite[pp.~119--120]{PR84} that the anticyclotomic pairing degenerates exactly when $M_0$ is symmetric.

\subsection{The quadratic form and its constant}
\begin{theorem}\label{thm:three}
Assume $\Sha(E/K)[p^\infty]$ finite.
\begin{enumerate}
\item There is $\tilde C_p\in W[1/p]$ with
\begin{equation}\label{eq:F2}
\mathcal F_2(a,b)=\tilde C_p\,\bigl[(a-b)^2\det H-(a+b)^2j^2\bigr]\qquad(a,b\in\Zp).
\end{equation}
\item If $\det H\ne0$, then
\begin{equation}\label{eq:Ctilde}
4\tilde C_p=\iop(f_0\omega_E)\,\Omega_p\,e^2\,(1-\alpha^{-1})^2\,q_p .
\end{equation}
\end{enumerate}
Thus the quadratic term of the Katz measure at $\psi$ is $q_p$ times the explicit quadratic form $\tfrac14\iop(f_0\omega_E)\Omega_pe^2(1-\alpha^{-1})^2\bigl[(a-b)^2\det H-(a+b)^2j^2\bigr]$. If $p$ is non-anomalous, $\tilde C_p\in W$.
\end{theorem}
\begin{proof}
(1) By the two-variable main conjecture \cite[Thm.~4.1(i)]{Rubin91}, Yager's theorem that the Katz measure generates the characteristic ideal of the local units modulo elliptic units \cite[Thm.~1]{Yager82}, and \cite[III.1.10--1.11]{deShalit}, the Katz measure generates the characteristic ideal of the dual of the Selmer group over $K_\infty$, twisted by $\kappa$ as in the proof of \cite[Thm.~12.2]{Rubin91}; hence $\int\kappa^{s}\kappa^{*s^*}d\mu_\psi=U\cdot\mathcal H_\fp(\kappa^s\kappa^{*s^*})$ with $U\in W[[S,T]]^\times$. Comparing homogeneous parts of degree two by Theorem~\ref{thm:PR} and Lemma~\ref{lem:dictionary}, and substituting $(s,s^*)=(ae,-be)$, gives \eqref{eq:F2}, with $\tilde C_p$ the product of $U(0)$, $\#\Sha(E/K)[\fp^\infty]$, the Euler factors of Theorem~\ref{thm:PR} and a power of $p$ from the substitution. None of this uses non-anomaly: when $p$ is anomalous an Euler factor is a non-unit and $\tilde C_p$ need not lie in $W$, but it lies in $W[1/p]$. The standing hypotheses of \cite{PR84} hold: the layers of $K_\infty/K$ are abelian over $K$, for which her Leopoldt-type hypotheses are Brumer's theorem.

(2) Since $\langle xy\rangle^{e}=(xy)^{e}$, one has $L^{\mathrm{Katz}}_p((1+p)^{et}-1)=\mathcal F(t,-t)$. Both sides vanish to order two, so comparing the coefficients of $t^2$ gives, exactly, $e^2\log_p(1+p)^2\,[T^2]L^{\mathrm{Katz}}_p=\mathcal F_2(1,-1)=4\tilde C_p\det H$ by \eqref{eq:F2}. By Lemma~\ref{lem:BK}, $[T^2]L^{\mathrm{Katz}}_p=f_0\omega_E\Omega_p\,c_2(p)$, and $c_2(p)=q_p(1-\alpha^{-1})^2\det H/\log_p(1+p)^2$; dividing by $\det H\ne0$ gives \eqref{eq:Ctilde}, for every $p$.
\end{proof}

The theorem does not use the value of $q_p$, and it holds whether or not the conjecture below is true: it says that the Katz measure, which is built from elliptic units and Hecke $L$-values, knows $q_p$. The rest of this section reads $q_p$ off it.

\begin{corollary}\label{prop:Cp}
Assume moreover that $p$ is non-anomalous and $\det H\ne0$, and put $X'=p^{-2}X\bmod p$ and $c=\iop(f_0\omega_E)(a_p-1)^2\in\F_p^\times$. Then
\begin{equation}\label{eq:three}
B'/2\equiv c\,q_p\det H',\qquad M'\equiv-c\,q_p\,j'^2,\qquad T_1'\equiv c\,q_p(\det H'-j'^2)\pmod p .
\end{equation}
In particular, the three congruences \eqref{eq:three} with $q_p$ replaced by $A_E$ are each equivalent to $q_p\equiv A_E\pmod p$, whenever the corresponding factor $\det H'$, $j'$ or $\det H'-j'^2$ is non-zero modulo $p$ and $p\nmid A_E$.
\end{corollary}
\begin{proof}
By Proposition~\ref{prop:moments}, $\mathcal F_2(2,0)\equiv\Omega_p^{2p-1}\iop(T_1)$ and $\mathcal F_2(1,1)\equiv\Omega_p^{2p-1}\iop(M)$ modulo $p^3W$. The functional equation of $L(\psib^{2p-1},s)$, whose root number is $+1$, gives $T_3=(\mathrm N\ff)^{p-1}T_1$; the edge value $T_1$ is moreover $\equiv4a_p\,\mathcal N(e,0)\pmod{\fp^3}$, because $G(u_1^2-1,0)-4G(u_1-1,0)\in p^3W$ for $G(S,0)\in S^2W[[S]]$ and $\Omega_p^{p-1}\equiv a_p^{-1}$. Since $\mathcal F_2(a,b)=\mathcal F_2(b,a)$, $B\equiv2(T_1-M)$ is $2\Omega_p^{-(2p-1)}\mathcal F_2(1,-1)$ modulo $p^3$. Put $C_p=4\Omega_p^{-(2p-1)}\tilde C_p\bmod\mathfrak m_W$; then \eqref{eq:F2} gives \eqref{eq:three} with $C_p$ in place of $c\,q_p$. Finally, by \eqref{eq:Ctilde}, $C_p\equiv\iop(f_0\omega_E)\Omega_p^{-(2p-2)}(1-\alpha^{-1})^2q_p$, and $\Omega_p^{-(2p-2)}(1-\alpha^{-1})^2\equiv a_p^2(1-a_p^{-1})^2=(a_p-1)^2$ by Lemma~\ref{lem:period}.
\end{proof}

The bracket $B$ is thus the parallelogram law of $\mathcal F_2$: the edge value $T_1$ is the $\fp$-direction, the central value $M$ the anticyclotomic one, and $B$ twice the cyclotomic one. The equivalence of the exact formula modulo $p$ with the cyclotomic congruence was found by Banwait \cite[Lemma~3.5, Cor.~6.9]{Banwait}; the point of Theorem~\ref{thm:three} is that the three directions carry the same constant, and that this constant is $q_p$ itself. The congruences \eqref{eq:three} also describe the primes where the criteria of \cite{CLS1} and \cite{Banwait} fail ($\det H'\equiv j'^2$ and $\det H'\equiv0$) and a third family ($j'\equiv0$); we do not pursue this here.

\subsection{Higher digits from central values}\label{ss:digits}
Corollary~\ref{prop:Cp} reads $q_p$ modulo $p$. The next digits seem to require the $p$-adic period $\Omega_p$ beyond its residue; they do not.
For $m\ge1$ let $c_k=(-1)^{k+1}\binom mk/k^2$ $(1\le k\le m)$; these solve the Vandermonde system $\sum_kc_kk^n=1$ for $n=2$ and $0$ for $3\le n\le m+1$, since $\sum_{k=0}^m(-1)^{m-k}\binom mkk^n=m!\,S(n,m)$ with $S(n,m)$ the Stirling numbers of the second kind.

\begin{theorem}\label{thm:digits}
Let $p$ be non-anomalous, $\Sha(E/K)[p^\infty]$ finite, $\det H\ne0$ and $v_p(j)=1$. Let $y_0\in\Zp$ be any lift of $a_p^{-2}\bmod p$ and $S_m=\sum_{k=1}^mc_k\,y_0^{\,k}\,\iop(M_k)$. Then
\begin{equation}\label{eq:digits}
q_p\equiv-\frac{S_m}{\iop(f_0\omega_E)\,(1-\alpha^{-1})^2\,e^2\,j^2}\pmod{p^{B_p(m)-2}},
\end{equation}
where $B_p(m)=v_p(m!)+\min_{n\ge m+2}\bigl(n-v_p(n!)\bigr)$.
For $p\ge m+3$ one has $B_p(m)=m+2$, so that $M_1,\dots,M_m$ determine $q_p$ modulo $p^m$; for $p=5$ and $m=2,4,7,12$ the precision is $5^{2},5^{3},5^{7},5^{12}$.
\end{theorem}
\begin{proof}
By Proposition~\ref{prop:moments}, $\mathcal F(k,k)=\Omega_p^{2ke+1}\iop(M_k)=\Omega_p\,y^k\iop(M_k)$ with $y=\Omega_p^{2e}\equiv y_0\pmod{pW}$ (Lemma~\ref{lem:period}). The function $g(t)=(y_0/y)^t\mathcal F(t,t)/(\Omega_pt^2)$ is analytic on $\Zp$, with $g(k)=y_0^k\iop(M_k)/k^2$ and Taylor coefficients $b_d$ of valuation $\ge d+2-v_p((d+2)!)$, because those of $\mathcal F(t,t)$ and of $(y_0/y)^t$ have valuations $\ge n-v_p(n!)$ and $\ge d-v_p(d!)$. Hence
\[S_m=\sum_k(-1)^{k+1}\tbinom mk g(k)=g(0)+(-1)^{m+1}m!\sum_{d\ge m}b_dS(d,m),\]
so $S_m\equiv g(0)=\Omega_p^{-1}\mathcal F_2(1,1)\pmod{p^{B_p(m)}}$. By \eqref{eq:F2} and \eqref{eq:Ctilde}, $\Omega_p^{-1}\mathcal F_2(1,1)=-4\Omega_p^{-1}\tilde C_pj^2=-\iop(f_0\omega_E)e^2(1-\alpha^{-1})^2q_p\,j^2$, of valuation $2+v_p(q_p)$.
\end{proof}

The right-hand side of \eqref{eq:digits} does not depend on the lift $y_0$ to the stated precision, which we checked numerically; a wrong residue for $y_0$ destroys the congruence. The values $M_k$ have weight $2ke+2$, and the cost grows with the weight (\S\ref{sec:evidence}).

\subsection{The conjecture}\label{ss:conj}
Let $\Reg_\infty$ be the N\'eron--Tate regulator of $P_1,P_2$ and
\[q_\infty=\frac{L''(E,1)}{2\,\Omega_E\,\Reg_\infty}\in\R .\]
Theorem~\ref{thm:three} and Theorem~\ref{thm:digits} make $q_p$ computable from critical Hecke $L$-values and $p$-adic heights, to any precision, for curves far beyond the reach of modular symbols, and the comparison with $q_\infty$ leads to the following statement.

\begin{conjecture}\label{conj:main}
Let $E=E_D$ have rank two and $p\equiv1\pmod4$ be a prime of good reduction with $\det H\ne0$. Then $q_\infty\in\Q$ and $q_p=q_\infty$ in $\Qp$.
\end{conjecture}

The complex Birch and Swinnerton-Dyer conjecture predicts $q_\infty=A_E$; with it, Conjecture~\ref{conj:main} is the exact formula $q_p=A_E$ of \cite{MTT}. Conversely, Conjecture~\ref{conj:main} is the exact formula of \cite{MTT} with $\#\Sha(E/\Q)$ replaced by its analytic value $\Shan=q_\infty\,\#E(\Q)_\tors^2/\Tam$. It involves only $L$-functions and heights, and it is of the same nature as the rank-one theorem of Rubin \cite[Thm.~8.2]{Rubin92}: for $E/\Q$ with complex multiplication and analytic rank one, and $y\in E(\Q)$ of infinite order,
\begin{equation}\label{eq:rubin-rank1}
(1-\alpha^{-1})^{-2}\,\frac{L_p'(\psi)}{\Omega_p\langle y,y\rangle_p}=\frac{L'(E,1)}{\Omega_\infty\langle y,y\rangle_\infty}\in K^\times,
\end{equation}
which says that the $p$-adic and complex quotients agree, without any arithmetic invariant. Exact Birch and Swinnerton-Dyer formulas with arithmetic invariants are known for certain CM curves of analytic rank zero \cite{BF24}, where the leading term is a rational number and all primes can be treated at once, and the $\fp$-part of the formula is known in analytic rank one for primes split in the CM field, under a condition on the conductor \cite[Thm.~A]{CastellaTNC}, and for curves over $\Q$ it follows from \eqref{eq:rubin-rank1} and the main conjecture \cite{Rubin91}; in rank two nothing exact is known. The main conjecture determines $q_p$ up to a unit because it is an equality of ideals; Sections~\ref{sec:rubin}--\ref{sec:AC} locate $q_p$ in two global classes and extract from identities of elements its square class.

\begin{remark}\label{rem:BKS}
Conjecture~\ref{conj:main} sits inside the Generalized Perrin-Riou Conjecture of Burns, Kurihara and Sano \cite[Conj.~4.8]{BKS}, which predicts the position of the Iwasawa-theoretic derivative $\kappa_\infty$ of Kato's zeta element. Their generalized Rubin formula \cite[Thm.~6.2]{BKS} determines the pairings of $\kappa_\infty$ with $E(\Q)$ in terms of $L_p(E,T)$, assuming only that $\Sha(E/\Q)[p^\infty]$ is finite, and by \cite[Cor.~6.7]{BKS} and its converse (valid when the $p$-adic regulator is non-zero) their conjecture for $E$ at $p$ is equivalent to the $p$-adic Beilinson formula $(1-\alpha^{-1})^{-1}(1-\alpha/p)L^{(2)}_{S,p}=L^*_S(E,1)R_p/(\Omega R_\infty)$, which is \eqref{eq:qp}$=q_\infty$ once the Euler factor $(1-\alpha^{-1})(1-\alpha/p)$ at $p$ is removed. Their formulation fixes an embedding $\R\hookrightarrow\C_p$ in order to compare $q_\infty$ with a $p$-adic number; this does not give the rationality of $q_\infty$, which remains a separate assertion.
\end{remark}

\begin{remark}\label{rem:units}
If $\Sha(E/\Q)$ is finite, the exact $p$-adic formula $q_p=A_E$ and the complex formula $q_\infty=A_E$ together imply Conjecture~\ref{conj:main}, and Conjecture~\ref{conj:main} together with either of them implies the other; Conjecture~\ref{conj:main} alone implies neither. None of these statements is easier than the others: in rank one the theorem \eqref{eq:rubin-rank1} compares the two analytic quotients directly, while the arithmetic forms are known only up to units, or prime by prime.
\end{remark}

\section{Elliptic units: the position of Rubin's class}\label{sec:rubin}

In rank $r\ge2$ Rubin \cite[\S\S4, 9, 11]{Rubin92} constructs from elliptic units classes $x^{(r)}_\fp\in E(K)\otimes\mathcal O_{K,\fp}$ and $x^{(r)}_\fpb$, and proves, exactly, that their formal logarithms and heights are given by derivatives of the Katz measure. We show that for $r=2$ these theorems determine the position of his class in $E(\Q)\otimes\Qp$ as $q_p$ times an explicit vector. Conjecture~\ref{conj:main} is then the statement that this position is the one given by $q_\infty$. The numerical test of the resulting identity is in \S\ref{ss:fourth-test}.

\subsection{Rubin's objects in our notation}
Take Rubin's field $F$ to be $K$, and $r=2$, in \cite{Rubin92}, and assume $\Sha(E/K)[p^\infty]$ finite. With $t,t^*$ as in Lemma~\ref{lem:dictionary}, the vectors $t(P_a)$ and $t^*(P_a)$ are bases of the modules $\check S=E(K)\otimes\mathcal O_{K,\fp}$ and $\check S^*$ of \cite[\S2]{Rubin92}. Rubin's pairing $\langle\ ,\ \rangle_\fp\colon\check S\times\check S^*\to\Qp$ \cite[Thm.~2.5]{Rubin92} is Perrin-Riou's \cite[IV, Lemme~19, Th\'eor\`eme~23]{PR84}, normalised by $\langle t(P),t^*(P)\rangle_\fp=h_\fp(P)$ \cite[III, Lemme~6]{PR84}, so by the proof of Lemma~\ref{lem:dictionary} its matrix is
\begin{equation}\label{eq:M}
\mathbf M=2M_0=\tfrac12\bigl(H+r^{-1}J\bigr),\qquad\det\mathbf M=\tfrac14\Delta,\qquad\Delta=\det H-j^2 .
\end{equation}
Let $\lambda$ be the formal logarithm of $\widehat E$ in $t=-x/y$, $\lambda_a=\lambda(P_a)$, and
\[\ell={}^t(\lambda_1,\lambda_2),\quad Q_\lambda=\lambda_2P_1-\lambda_1P_2\in E(\Q)\otimes\Qp,\quad h^+(Q_\lambda)={}^t\ell\,\adj(H)\,\ell .\]
$Q_\lambda$ spans the line of vanishing logarithm, and $\log_{E,\fp}(t(P_a))=\log_{E,\fpb}(t^*(P_a))=\lambda_a$ because $\lambda$ has rational coefficients. Let $x^*=x^{(2)}_\fpb\in\check S^*$ be the class of \cite[\S\S6, 11]{Rubin92} with $\fp$ and $\fpb$ interchanged \cite[\S10]{Rubin92}; it is non-zero when $\Delta\ne0$ \cite[\S11]{Rubin92}. Finally let
\[G^*(a,b)=\int x^{ae}y^{1-be}\,d\mu=\Omega_p^{(a+b)e-1}\,\iop\mathcal N(ae-1,be-1)\quad(a\ge b\ge1),\]
the branch of the Katz measure through the conjugate character $\psi^*=\psib$, extended to an Iwasawa function of $(a,b)\in\Zp^2$; it vanishes at the origin, and we write $\kappa_1(a-b)$ for the linear part of $\Omega_pG^*$.

\begin{lemma}\label{lem:rubin-norm}
The two-variable $p$-adic $L$-function $\Lc_\fp$ of \cite[Thm.~7.1]{Rubin92} is the measure $\mu$ of Proposition~\ref{prop:moments}, and its derivatives at $\psi$ and $\psi^*$ along the $\fp$-direction are
\begin{equation}\label{eq:L2}
\Lc_\fp^{(2)}(\psi)=\frac{\mathcal F_2(1,0)}{e^2}=\tfrac14\,q_p\,(1-\alpha^{-1})^2\,\iop(f_0\omega_E)\,\Omega_p\,\Delta,\qquad\Omega_p\,\Lc_\fpb^{(1)}(\psi^*)=\frac{\kappa_1}{e}.
\end{equation}
\end{lemma}
\begin{proof}
Rubin's interpolation formula at $\psi^k\psi^{*j}$, $0\le-j<k$, coincides term by term with \eqref{eq:EK}--\eqref{eq:moments}: his Euler factor is $1-\Frob^{w}/p^{l+1}$, his $L$-value has the factor at $\fpb$ removed, $\sqrt{-D_K}/2\pi=1/\pi$, his lattice $\Omega\ff$ is $\Omega_\infty\OK$, and his $\Omega_p$ is the inverse period. The first formula follows from $\mathcal F_2(1,0)=\tilde C_p\Delta$ \eqref{eq:F2} and \eqref{eq:Ctilde}; the second is the definition of $\kappa_1$. Rubin's derivatives are Taylor coefficients in the variable $s$ of $\psi\langle\psi\rangle^s$, normalised by $\log_p\psi(\gamma)=p$; another normalisation would change \eqref{eq:fourth} below by a power of $p$, which the data of \S\ref{ss:fourth-test} exclude.
\end{proof}

\begin{theorem}[Position of Rubin's class]\label{thm:position}
Let $p\equiv1\pmod4$ be of good reduction $($possibly anomalous$)$, $\Sha(E/K)[p^\infty]$ finite, $\det H\ne0$ and $\Delta\ne0$. Then:
\begin{enumerate}
\item[(i)] for every $y\in\check S$, $\displaystyle\langle y,x^*\rangle_\fp=\frac{\iop(f_0\omega_E)}{4p}\Bigl(1-\frac\alpha p\Bigr)\,q_p\,\Delta\,\log_{E,\fp}(y)$;
\item[(ii)] in the basis $t^*(P_1),t^*(P_2)$,
\begin{equation}\label{eq:position}
x^*=\frac{\iop(f_0\omega_E)}{2p}\Bigl(1-\frac\alpha p\Bigr)\,q_p\,\adj\bigl(H+r^{-1}J\bigr)\,\ell ;
\end{equation}
\item[(iii)] $\log_{E,\fpb}(x^*)=\frac{\iop(f_0\omega_E)}{2p}(1-\frac\alpha p)\,q_p\,h^+(Q_\lambda)$, and $\log_{E,\fp}(x^{(2)}_\fp)=\frac{\mathrm N\ff}{p\,e}\bigl(1-\frac\alpha p\bigr)^{-1}\kappa_1$.
\end{enumerate}
\end{theorem}
\begin{proof}
For anomalous $p$ we use \eqref{eq:F2} and \eqref{eq:Ctilde}, which hold for every $p\equiv1\pmod4$ of good reduction with $\tilde C_p\in W[1/p]$ (Theorem~\ref{thm:three}); non-anomaly is only needed for the reductions modulo $p$ in Corollary~\ref{prop:Cp}.
(i) Write $x^*=\Phi^*(\xi^*)$ as in \cite[Prop.~2.4, \S6]{Rubin92} and let $\alpha_w$ be the lift of \cite[Prop.~9.4]{Rubin92}. With Rubin's field equal to $K$, Theorem~2.5 and Proposition~5.1 of \cite{Rubin92} give $p^{-1}\langle y,x^*\rangle_\fp=\log_{E,\fp}(y)\,(1-\alpha/p)\,\delta_w(\alpha_w)$, and \cite[Prop.~9.4(i)]{Rubin92} gives $\delta_w(\alpha_w)=p^{-2}(1-\alpha^{-1})^{-2}\Omega_p^{-1}\Lc^{(2)}_\fp(\psi)$; substitute \eqref{eq:L2}. (The printed statement of \cite[Thm.~9.3(ii)]{Rubin92} carries a factor $p(p-1)$ where Theorem~2.5 and Proposition~5.1 give $p$; the latter is also what \cite[Thm.~9.3(i), Cor.~9.2]{Rubin92} give, and it is confirmed numerically in \S\ref{ss:fourth-test}.)
(ii) Applying (i) to $y=t(P_1),t(P_2)$ gives $\mathbf MX^*=\tfrac14\iop(f_0\omega_E)p^{-1}(1-\alpha/p)q_p\Delta\,\ell$ for the coordinates $X^*$ of $x^*$, and $\mathbf M^{-1}=2\adj(H+r^{-1}J)/\Delta$.
(iii) The first formula is ${}^t\ell X^*$, since ${}^t\ell\adj(H+r^{-1}J)\ell=h^+(Q_\lambda)$ by antisymmetry of $\adj(J)$; the second is \cite[Thm.~9.5(ii)]{Rubin92} with $r=2$ and trivial character, $\log_{E,\fp}(x^{(2)}_\fp)=p^{-1}\mathrm N\ff\,(1-\psi^*(\fp)/p)^{-1}\Omega_p\Lc^{(1)}_\fpb(\psi^*)$, with $\iop\psi^*(\fp)=\alpha$ and \eqref{eq:L2}.
\end{proof}

The two logarithms in (iii) are linked by a fourth root of unity. Both classes are built from one orbit $B$ of \cite[p.~335]{Rubin92}, and complex conjugation $\tau$ interchanges every ingredient of the two constructions except the orbit, which it sends to $u_EB$ with $u_E=\bar\Omega/\Omega$, $\Omega=\Omega_\infty/f_0$ \cite[Prop.~6.1]{Rubin92}. By \cite[Prop.~3.1(iv) and the Remark on p.~336]{Rubin92} and the $\Zp$-linearity of the constructions \cite[Lemma~4.1, Thm.~4.2, Prop.~2.4]{Rubin92},
\begin{equation}\label{eq:tau}
\tau(x^*)=\iop(u_E)\,x^{(2)}_\fp,\qquad\text{so}\qquad\log_{E,\fpb}(x^*)=\iop(u_E)\log_{E,\fp}(x^{(2)}_\fp).
\end{equation}
The same constant is obtained by comparing the interpolation formulas of $\Lc_\fp$ and of its conjugate $\Lc_\fpb$, which differ by $(\Omega/\bar\Omega)^{j-k}$ and the complex conjugation of all the other factors. Explicitly $\iop(u_E)=\iop(i^{e_f})$ if $D>0$ and $\iop(i^{e_f-1})$ if $D<0$; it is $-1$ for $E_{56}$ and $-i$ for $E_{17}$.

\begin{theorem}\label{thm:fourth}
Under the hypotheses of Theorem~\ref{thm:position},
\begin{equation}\label{eq:fourth}
\frac{\kappa_1}{p-1}=\frac12\,\iop\Bigl(\frac{\bar\omega_E}{f_0}\Bigr)\Bigl(1-\frac\alpha p\Bigr)^2\,q_p\,h^+(Q_\lambda).
\end{equation}
\end{theorem}
\begin{proof}
Equate the two expressions of Theorem~\ref{thm:position}(iii) through \eqref{eq:tau}, and use $f_0\omega_E\,\Omega/\bar\Omega=\bar f_0\bar\omega_E$ (because $\omega_E\Omega_\infty=\Omega_E$ is real) and $f_0\bar f_0=\mathrm N\ff$.
\end{proof}

\begin{corollary}\label{cor:position}
Under the hypotheses of Theorem~\ref{thm:position}, Conjecture~\ref{conj:main} at $p$ is equivalent to: $q_\infty\in\Q$ and \eqref{eq:position} holds with $\iop(q_\infty)$ in place of $q_p$; equivalently, $\langle t(P),x^*\rangle_\fp=\tfrac14\iop(f_0\omega_E)p^{-1}(1-\alpha/p)\Delta\,\lambda(P)\,\iop(q_\infty)$ for one rational point $P$ of infinite order. If $h^+(Q_\lambda)\ne0$ it is also equivalent to \eqref{eq:fourth} with $q_\infty$ in place of $q_p$.
\end{corollary}

This is the rank-two analogue of Rubin's rank-one formula \cite[Cor.~10.3, Thm.~10.4]{Rubin92}, which computes a rational point from the $p$-adic $L$-function; Rubin himself observes \cite[Remark~(1), p.~348]{Rubin92} that the $p$-adic Birch and Swinnerton-Dyer conjecture predicts the value of the quotient $(1-\alpha^{-1})^{-2}\Lc^{(2)}_\fp(\psi)/(\tilde\Omega_\fp R_\fp)$, and by \eqref{eq:L2} that quotient is an explicit multiple of $q_p$. Castella \cite{Castella26} shows that the generalised Kato class of a CM curve of rank two, when non-zero, spans the line of $Q_\lambda$, and that, if $\Sha(E/\Q)[p^\infty]$ is finite, $L_p$ vanishes to order two and the $\psi^*$-branch to order exactly one (i.e.\ $h^+(Q_\lambda)\ne0$), it spans the Selmer group together with Rubin's class; Theorem~\ref{thm:position} gives the second vector of his basis explicitly.

\section{Derived Heegner classes over an auxiliary field}\label{sec:heegner}

Over $K=\Q(i)$ there are no Heegner points for $E$: the Heegner hypothesis fails at $2$, since $\Z[i]/(1+i)^k$ is not cyclic for $k\ge2$ while $2^5\mid N_E$. We therefore pass to an auxiliary imaginary quadratic field $F$ over which $E$ has rank three and Heegner points exist, and use Howard's $\Lambda$-adic Gross--Zagier theorem to express $q_p$ through the derivative of the Heegner class.

\subsection{Hypotheses}\label{ss:hyp}
Let $F=\Q(\sqrt{-d})$, $\chi=\chi_F$, $E^\chi$ the quadratic twist, $\mathcal K=KF=\Q(i,\sqrt{-d})$ and $F^+=\Q(\sqrt d)$.

\begin{hypothesis}\label{hyp}
\begin{enumerate}
\item[(H1)] $E(\Q)$ has rank two and $\ord_{s=1}L(E,s)=2$; $\ord_{s=1}L(E^\chi,s)=1$.
\item[(H2)] $p\equiv1\pmod4$ and $p\nmid N_E$.
\item[(H3)] $-d$ is an odd fundamental discriminant; every prime dividing $N_E$ splits in $F$, and $p$ splits in $F$.
\item[(H4)] $p\nmid6\,h^-_{\mathcal K}\,D_{F^+}$, where $h^-_{\mathcal K}=h_{\mathcal K}/h_{F^+}$.
\item[(H5)] $c_2(p)\ne0$.
\item[(H6)] $Q_2\ne0$, with $Q_2$ as in \eqref{eq:Q2} below.
\end{enumerate}
\end{hypothesis}

All of them are checked by a finite computation; for (H5) one needs $c_2(p)$ modulo $p^{v_p(c_2(p))+1}$. The sign of $L(E/F,s)=L(E,s)L(E^\chi,s)$ is $-1$ by (H3), consistently with (H1).

\begin{lemma}\label{lem:consequences}
Under Hypotheses~\ref{hyp}:
\begin{enumerate}
\item $E^\chi(\Q)$ has rank one, generated modulo torsion by a point $P_3$ that we view in $E(F)$, and $E(F)\otimes\Q=(E(\Q)\oplus E^\chi(\Q))\otimes\Q$ has dimension three; $E(F)[p]=0$.
\item $\Sha(E/\Q)[p^\infty]$, $\Sha(E^\chi/\Q)[p^\infty]$ and $\Sha(E/F)[p^\infty]$ are finite, and $\det H\ne0$.
\item The cyclotomic $p$-adic height $h_3$ of $P_3$ is non-zero.
\item The residual representation $\bar\rho\colon G_F\to\mathrm{Aut}(E[p])$ is irreducible and $p$-distinguished.
\end{enumerate}
\end{lemma}
\begin{proof}
(1) is Gross--Zagier \cite{GZ86} and Kolyvagin \cite[Cor.~C]{Kolyvagin90} for $E^\chi$, and $E(F)[p]=0$ because $E(F)_\tors$ is a $2$-group for $p\ge5$. (2): by (H5), Rubin's cyclotomic main conjecture \cite[Thm.~12.3]{Rubin91} and the control theorem, the corank of $\mathrm{Sel}_{p^\infty}(E/\Q)$ is at most two, and Perrin-Riou's theorem \cite[Thm.~V.8]{PR84} gives $\det H\ne0$ (\S\ref{sec:constant}); for $E^\chi$ use Kolyvagin; and $[F:\Q]=2$ is prime to $p$. (3): $h_3$ is twice the $\fp$-adic height of $P_3$ for $K$, because the anticyclotomic height of a point defined over $\Q$ vanishes; the latter is non-zero by Bertrand's theorem \cite[Cor.~1]{Bertrand84}, in the form \cite[Appendix, Thm.~A.1]{AH06}, as $E^\chi$ has complex multiplication by $\OK$. (4): $\bar\rho|_{G_{\mathcal K}}=\bar\chi_\fp\oplus\bar\chi_\fpb$, and $\bar\rho|_{G_F}$ is induced from $\bar\chi_\fp$; by Mackey's criterion it is irreducible if $\bar\chi_\fp\ne\bar\chi_\fpb$ on $G_{\mathcal K}$. As $p\nmid d$, the extension $\mathcal K/K$ is unramified above $\fp$, and the inertia group at $\fp$ acts on $E[\fp]$ through the cyclotomic character $\omega\ne1$ and trivially on $E[\fpb]$. The same computation gives the distinction at $p$.
\end{proof}

We checked (4) also numerically: for $(E,F,p)=(E_{56},\Q(\sqrt{-31}),5)$, $(E_{-289},\Q(\sqrt{-111}),5)$, $(E_{-68},\Q(\sqrt{-111}),5)$ the $5$-division polynomial factors over $F$ into irreducible factors of degrees $4$ and $8$, and for $(E_{56},\Q(\sqrt{-55}),13)$ the $13$-division polynomial into factors of degrees $12$ and $72$.

\subsection{\texorpdfstring{Heights over $F$}{Heights over F}}
Let $F_\infty$ be the $\Zp^2$-extension of $F$, with cyclotomic and anticyclotomic $\Zp$-extensions $F^{c}_\infty$, $F^{a}_\infty$; let $\gamma_0$ be the generator of $\mathrm{Gal}(F^c_\infty/F)$ with $\chi_{\rm cyc}(\gamma_0)=1+p$, $\sigma_a$ a topological generator of $\mathrm{Gal}(F^a_\infty/F)$, $\Lambda_F=\Zp[[T_c,T_a]]$ with $T_c=\gamma_0-1$, $T_a=\sigma_a-1$, and $\ell_a\colon\mathrm{Gal}(F^a_\infty/F)\to\Qp$ a non-zero homomorphism. Let $h^F_c$ and $h^F_a$ be the cyclotomic and anticyclotomic $p$-adic height pairings on $E(F)$, for the homomorphisms $\log_p\circ\chi_{\rm cyc}$ and $\ell_a$, normalised over $F$, so that $h^F_c=2h^+$ on $E(\Q)$. Complex conjugation changes the sign of $\ell_a$ and fixes $E(\Q)$, and multiplies $E^\chi(\Q)$ by $-1$; hence $h^F_a$ vanishes on $E(\Q)\times E(\Q)$ and on $E^\chi(\Q)\times E^\chi(\Q)$, and $h^F_c$ vanishes on $E(\Q)\times E^\chi(\Q)$. In the basis $P_1,P_2,P_3$ of $E(F)\otimes\Q$,
\[H_c=\begin{pmatrix}2H&0\\0&2h_3\end{pmatrix},\qquad H_a=\begin{pmatrix}0&v\\{}^tv&0\end{pmatrix},\qquad v=\begin{pmatrix}h^F_a(P_1,P_3)\\h^F_a(P_2,P_3)\end{pmatrix},\]
and
\begin{equation}\label{eq:Q2}
\det\bigl(XH_c+YH_a\bigr)=X^3\,8h_3\det H+XY^2\,Q_2,\qquad Q_2=-{}^tv\,\adj(2H)\,v .
\end{equation}
Put $Q_v=v_2P_1-v_1P_2\in E(\Q)\otimes\Qp$; then $h^F_a(Q_v,P_3)=0$, $h^F_a(Q_v,E(\Q))=0$ and $h^F_c(Q_v,Q_v)={}^tv\adj(2H)v=-Q_2$.

\subsection{\texorpdfstring{The main conjecture on the plane of $F$}{The main conjecture on the plane of F}}
Let $X$ be the Pontryagin dual of the $p^\infty$-Selmer group of $E$ over $F_\infty$, a $\Lambda_F$-module, and $L_{\mathcal K}\in\Lambda_F\otimes W$ the restriction to the plane of $F$ of the Katz $p$-adic $L$-function of $\mathcal K$ attached to the Hecke character of $E$ over $\mathcal K$, as constructed by Hida and Tilouine \cite[Thm.~II]{HT93} and used by Hsieh \cite[\S6.5]{Hsieh14}.

\begin{lemma}\label{lem:euler}
Let $A=\mathcal O[[T_1,\dots,T_r]]$, $f\in A$ with $A/f\simeq\mathcal O[[T_1,\dots,T_{r-1}]]$, and $X$ a finitely generated torsion $A$-module. Then $X/fX$ is $A/f$-torsion if and only if $f\nmid\mathrm{char}\,X$, and in that case $\mathrm{char}_{A/f}(X/fX)=(\mathrm{char}_AX\bmod f)\cdot\mathrm{char}_{A/f}(X[f])$.
\end{lemma}
\begin{proof}
Both sides are determined by their orders at the height-one primes of $A/f$, i.e.\ at the height-two primes $\mathfrak q\supset(f)$ of $A$. Localise at such a $\mathfrak q$: if $f\nmid\mathrm{char}\,X$, then $X_\mathfrak q/f$ and $X_\mathfrak q[f]$ have finite length, and $\chi_f(N)=\ell(N/fN)-\ell(N[f])$ is additive and vanishes on modules of finite length, so it may be computed on an elementary module $\bigoplus A/(g_i)$, where it is $\ord_{\mathfrak q}(\prod g_i\bmod f)$. If $f\mid\mathrm{char}\,X$, then $X$ has a quotient pseudo-isomorphic to $A/(f)$, and $X/fX$ is not torsion.
\end{proof}

\begin{theorem}\label{thm:planeMC}
Let $R$ be the localisation of $\Lambda_F\otimes W[1/p]$ at the augmentation ideal. Under Hypotheses~\ref{hyp}, $X$ is $\Lambda_F$-torsion and $\mathrm{char}_{\Lambda_F}(X)\,R=(L_{\mathcal K})\,R$.
\end{theorem}
\begin{proof}
We separate (a) the control from three variables to the plane, (b) the comparison of local conditions, and (c) the control from the plane to the cyclotomic line. In each step what is needed is that the error terms have characteristic ideals prime to the augmentation ideal, i.e.\ that they are units in $R$; we do not claim that they vanish.

\emph{Torsion.} Let $X_{\rm cyc}$ be the dual Selmer group of $E$ over $F^c_\infty$, that is, of $E$ and $E^\chi$ over $\Q_\infty$. By \cite[Thm.~12.3]{Rubin91}, $\mathrm{char}(X_{\rm cyc})=(L_p(E)L_p(E^\chi))$, which is non-zero by (H5) and Lemma~\ref{lem:consequences}(3). By (c) below, $X/T_aX$ and $X_{\rm cyc}$ differ by modules that are finite, so $X/T_aX$ is torsion over $\Lambda_F/(T_a)$, i.e.\ $(X/T_aX)\otimes\mathrm{Frac}(\Lambda_F/(T_a))=0$. This says $X_{(T_a)}/T_aX_{(T_a)}=0$, so $X_{(T_a)}=0$ by Nakayama's lemma over the local ring $(\Lambda_F)_{(T_a)}$. Hence $X\otimes\mathrm{Frac}(\Lambda_F)=0$, i.e.\ $X$ is torsion, and $T_a\nmid\mathrm{char}\,X$.

\emph{Divisibility.} Let $X_3$ be the dual of the Selmer group of the Hecke character of $E$ over the $\Zp^3$-extension of $\mathcal K$, as in \cite{Hsieh14}. By \cite[Thm.~8.16, Cor.~8.19]{Hsieh14}, under (H2) and (H4), $\mathrm{char}(X_3)\subset(L_3)$, with $L_3$ the Katz $p$-adic $L$-function. (a) The plane of $F$ is the quotient of the $\Zp^3$-extension of $\mathcal K$ by the $\chi_{-4}$-part of its Galois group, which lies in the minus part; this is case (ii) of \cite[Remark~8.4]{Hsieh14}, with $d=[F^+:\Q]=2$ and $t=1$. Let $f_3\in\mathcal O[[T_1,T_2,T_3]]$ be the corresponding generator. By \cite[Lemma~8.3]{Hsieh14}, $X_3/f_3X_3$ and the Selmer group of the plane with the local conditions of \cite{Hsieh14} differ by local terms whose characteristic ideals are prime to the augmentation ideal: $\bar\chi_\fp=\omega$ on inertia at $\fp$, so the global term is a unit; the inertia at the bad primes acts through a group of order prime to $p$; and $\alpha$ is not a root of unity. (b) The local conditions of \cite{Hsieh14} and the ordinary conditions defining $X$ differ at the places above $p$ by terms whose specialisations at the augmentation ideal are $H^0$ and $H^2$ of $F_v$ with coefficients in the unramified quotient $V^-$ or its dual, on which Frobenius acts by $\alpha^{\pm1}\ne1$; so these specialisations vanish and the characteristic ideals are prime to the augmentation ideal. As $X$ is torsion, so is $X_3/f_3X_3$, hence $f_3\nmid\mathrm{char}(X_3)$ by Lemma~\ref{lem:euler} ($r=3$), and the lemma gives $\mathrm{char}(X_3/f_3X_3)\subset(\mathrm{char}(X_3)\bmod f_3)\subset(L_{\mathcal K})$, the term $\mathrm{char}(X_3[f_3])$ only adding a further factor. With (a) and (b), $\mathrm{char}(X)R\subset(L_{\mathcal K})R$.

\emph{Equality.} Write $\mathrm{char}(X)R=(L_{\mathcal K}\,G)R$ with $G\in R$. (c) The restriction to the cyclotomic line $T_a=0$ is controlled by \cite[Lemma~8.3]{Hsieh14}; there the local terms are finite, since a non-trivial cyclotomic character is ramified at $p$. As $T_a\nmid\mathrm{char}\,X$, Lemma~\ref{lem:euler} ($r=2$) gives
\[\mathrm{char}(X_{\rm cyc})=(L_{\mathcal K}|_{\rm cyc})\cdot(G\bmod T_a)\cdot\mathrm{char}(X[T_a])\]
after localisation at the augmentation ideal of the cyclotomic Iwasawa algebra. Now $(L_{\mathcal K}|_{\rm cyc})=(L_p(E)L_p(E^\chi))$ by the factorisation in the proof of \cite[Thm.~8.18]{Hsieh14} together with \cite[p.~67]{Rubin91}, and this equals $\mathrm{char}(X_{\rm cyc})$, which is non-zero. Cancelling it, $(G\bmod T_a)\cdot\mathrm{char}(X[T_a])$ is a unit there. Hence $G\bmod T_a$ is a unit at the origin, $G(0)\ne0$, and $G$ is a unit in $R$. (The non-vanishing of $L_{\mathcal K}|_{\rm cyc}$ alone would not suffice: it is the equality of the two cyclotomic characteristic ideals that forces the extra factor to be a unit.)
\end{proof}

\subsection{The cubic form}
Let $\mathfrak m=(T_c,T_a)$.

\begin{theorem}\label{thm:cubic}
Under Hypotheses~\ref{hyp} there is $C\in W[1/p]^\times$ such that
\[L_{\mathcal K}(T_c,T_a)\equiv C\,\det\Bigl(\frac{T_c}{\log_p(1+p)}\,H_c+\frac{T_a}{\ell_a(\sigma_a)}\,H_a\Bigr)\pmod{\mathfrak m^4}.\]
In particular $[T_c^2T_a]L_{\mathcal K}=[T_a^3]L_{\mathcal K}=0$ and
\begin{equation}\label{eq:ratio}
\frac{[T_cT_a^2]\,L_{\mathcal K}}{[T_c^3]\,L_{\mathcal K}}=\frac{\log_p(1+p)^2\,Q_2}{\ell_a(\sigma_a)^2\,8h_3\det H}.
\end{equation}
\end{theorem}
\begin{proof}
Let $C^\bullet=\widetilde{R\Gamma}_{f,\mathrm{Iw}}(F_\infty/F,T_pE)$ be Nekov\'a\v{r}'s Selmer complex with ordinary local conditions, and $R$ the localisation of $\Lambda_F\otimes W[1/p]$ at $\mathfrak m$, as in Theorem~\ref{thm:planeMC}. By the control theorem \cite[8.10.1]{Nekovar06}, $C^\bullet\otimes^{\mathbf L}_R R/\mathfrak m$ computes $H^1_f(F,V_pE)=E(F)\otimes\Qp$ in degree one, its dual in degree two, and nothing else, by Lemma~\ref{lem:consequences}(1)--(2). So $C^\bullet_R$ is quasi-isomorphic to a minimal complex $[R^3\xrightarrow{\Phi}R^3]$ with $\Phi\in\mathfrak m\,M_3(R)$, and $\mathrm{char}(X)R=(\det\Phi)$; the error terms at the bad primes are units at the origin because the Tamagawa numbers of $E$ over $F$ are powers of $2$ and the bad primes do not ramify in $F_\infty$. Write $\Phi\equiv T_cB_c+T_aB_a\pmod{\mathfrak m^2}$. The Bockstein map $H^1\to H^2\otimes\mathfrak m/\mathfrak m^2$ is $\Phi\bmod\mathfrak m^2$, and by \cite[11.1.4, 11.3.9]{Nekovar06} the associated pairings are, up to the sign $h^{\rm norm}=-\tilde h_{\pi,1,1}$ of \cite[11.3.9]{Nekovar06}, the $p$-adic heights for the homomorphisms dual to $T_c$ and $T_a$, i.e.\ $h^F_c/\log_p(1+p)$ and $h^F_a/\ell_a(\sigma_a)$; these coincide with the classical heights of Schneider and Mazur--Tate by \cite[8.1]{Nekovar93}. Hence $B_c=H_cG^{-1}/\log_p(1+p)$, $B_a=H_aG^{-1}/\ell_a(\sigma_a)$ for a fixed matrix $G$ of the duality, and $\det\Phi\equiv\det(G)^{-1}\det\bigl(T_cH_c/\log_p(1+p)+T_aH_a/\ell_a(\sigma_a)\bigr)\pmod{\mathfrak m^4}$. Theorem~\ref{thm:planeMC} transfers this to $L_{\mathcal K}$, the unit only changing $C$; finally \eqref{eq:Q2}. A common sign in the comparison of the Bockstein pairings with the heights (even different signs in the cyclotomic and anticyclotomic directions) only changes $C$ by a sign, since \eqref{eq:Q2} is even in $Y$; so it does not affect \eqref{eq:ratio}.
\end{proof}

Unlike a descent along each line separately, the minimal complex gives a single constant for all directions, which is what makes \eqref{eq:ratio} exact.

\subsection{The Hida--Perrin-Riou function}
Let $f$ be the newform of $E$, and $L_f\in\Zp[[\Gamma]]\otimes\Qp$, $\Gamma=\mathrm{Gal}(F_\infty/F)$, the two-variable $p$-adic $L$-function of Hida and Perrin-Riou \cite[Thm.~1.1]{PR87}, restricted to the plane, with the normalisation of \cite{Howard05}. At a finite-order character $\mathcal W$ of $\Gamma$, evaluated on ideals through the Artin map, with $\mathrm N\mathfrak f(\mathcal W)=p^{n_1+n_2}$, it interpolates \cite[pp.~457--458]{PR87}
\[\bar{\mathcal W}(\mathfrak D')\,\frac{\sqrt{|D_F|}\,p^{(n_1+n_2)/2}}{\alpha^{n_1+n_2}}\,V_p(\mathcal W)\,w(\mathcal W)\,\frac{L(f,\bar{\mathcal W},1)}{\Omega_f},\qquad\Omega_f=8\pi^2\int|f|^2dx\,dy,\]
with $V_p(\mathcal W)$ an Euler factor at $p$ and $w(\mathcal W)$ the Artin root number. Let $\vartheta=\psi\circ\mathrm N_{\mathcal K/K}$, $\widetilde{\mathcal W}=\mathcal W\circ\mathrm N_{\mathcal K/F}$ and $\lambda_{\mathcal W}=\vartheta^{-1}\widetilde{\mathcal W}^{-1}$, of infinity type $\Sigma$ in the sense of \cite[(0.3)]{HT93}, $\Sigma$ being induced from $K$. Let $\mu_{\rm HT}$ be the Katz measure of \cite[Thm.~II]{HT93} with auxiliary conductor $\mathfrak C$ the prime-to-$p$ conductor of $\vartheta$, and put
\[L^{\rm PR}_{\mathcal K}(\mathcal W)=\int\widehat{\lambda_{\mathcal W}}\,d\mu_{\rm HT} .\]

\begin{proposition}\label{prop:LfLK}
There are $c_0\in W[1/p]^\times$ and $\sigma_0\in\Gamma$ such that
\begin{equation}\label{eq:LfLK}
L_f=c_0\,[\sigma_0]\,L^{\rm PR}_{\mathcal K}.
\end{equation}
\end{proposition}
\begin{proof}
\emph{Complex values.} Since $\vartheta\vartheta^a=\mathrm N$ for the non-trivial $a\in\mathrm{Gal}(\mathcal K/F)$ and $\widetilde{\mathcal W}^a=\widetilde{\mathcal W}$, $L(0,\lambda_{\mathcal W})=L(1,\vartheta\widetilde{\mathcal W}^{-1})=L(f,\bar{\mathcal W},1)$.

\emph{Factors at $p$.} Let $\mathcal W$ be ramified at both primes $v_1,v_2$ of $F$ above $p$; then $V_p(\mathcal W)=1$. At each $v_j$, with $\mu_j=\mathcal W_{v_j}$ of conductor $p^{n_j}$ and $w_j$ the place of $\Sigma$ above $v_j$, the local Gauss sum of \cite[(0.8), (0.10)]{HT93} is
\[p^{-n_j}G(2\delta;\lambda_{\mathcal W,w_j})=\alpha^{-n_j}\mu_j(p)^{n_j}\mu_j((2\delta)_{w_j})^{-1}\sum_{u\in(\Z/p^{n_j})^\times}\mu_j(u)^{-1}e_p(u/p^{n_j}),\]
by $\vartheta_{w_j}(p)=\beta=p/\alpha$ and the substitution $u\mapsto(2\delta)_{w_j}u$, $(2\delta)_{w_j}\in\Zp^\times$. This is $\alpha^{-n_j}p^{n_j/2}$ times the local root number of $\mathcal W$ at $v_j$, up to $\mu_j((2\delta)_{w_j})^{-1}$; so the factors at $p$ agree up to $\mathcal W$ at a fixed id\`ele.

\emph{Other factors.} Since $\widetilde{\mathcal W}$ is unramified outside $p$, every $\lambda_{\mathcal W}$ is primitive of conductor $\mathfrak C$ away from $p$, and the Euler factors at the primes dividing $\mathfrak C$ in \cite[Thm.~II]{HT93} equal $1$, with the convention $\lambda(\mathfrak L)=0$ for $\mathfrak L\mid\mathfrak C$ of loc.\ cit. The remaining factors away from $p$ (root numbers, the different, the change of additive character) are constant or equal to $\mathcal W$ at a fixed id\`ele. The infinity type is fixed, so the periods and the $\Gamma$-factors are constant. The ratio of the complex periods of the two formulas is algebraic, as one sees by comparing them at one character $\mathcal W$, ramified at $v_1$ and $v_2$, with $L(f,\bar{\mathcal W},1)\ne0$; it exists because $L_f\ne0$ (by \eqref{eq:T3} below) and by the uniqueness statement below. The $p$-adic period of \cite{HT93} remains a fixed factor in $W^\times$.

\emph{Uniqueness.} A bounded measure on $\Gamma$ that vanishes at every finite-order character ramified at both $v_1$ and $v_2$ is zero: multiplied by $([g_1]-1)([g_2]-1)$, with $g_j\ne1$ in the inertia group at $v_j$, it vanishes at all finite-order characters, and the Iwasawa algebra of $\Zp^2$ is a domain. Applied to $L_f-c_0[\sigma_0]L^{\rm PR}_{\mathcal K}$, this gives \eqref{eq:LfLK} for some $c_0\ne0$ and $\sigma_0\in\Gamma$; and $c_0\in W[1/p]^\times$, since both sides lie in $\Zp[[\Gamma]]\otimes W[1/p]$ and have a non-zero coefficient in common.
\end{proof}

The branch of the Katz measure in Proposition~\ref{prop:LfLK} is the one of \cite[\S6.5]{Hsieh14}. With the geometric reciprocity map of \cite[\S1.3]{Hsieh14}, the character $\rho_\fp$ of $G_{\mathcal K}$ on $T_\fp E$ satisfies $\rho_\fp(\mathrm{Frob}^{\rm geom}_w)=\iota_p(\vartheta(w))^{-1}=\widehat{\vartheta^{-1}}(\mathrm{Frob}^{\rm geom}_w)$ for $w$ outside the conductors, and Chebotarev identifies the two characters. So $L^{\rm PR}_{\mathcal K}$ and $L_{\mathcal K}$ are the same branch, centred at $\widehat{\vartheta^{-1}}$, and can differ only in the parametrisation of $\Gamma$.
Depending on whether Galois characters are identified with id\`ele class characters by the arithmetic or the geometric reciprocity map, $L^{\rm PR}_{\mathcal K}$ is $L_{\mathcal K}$ or $L_{\mathcal K}^\iota$, where $\iota(\gamma)=\gamma^{-1}$ acts by $(T_c,T_a)\mapsto(-T_c,-T_a)$ to first order. Either way \eqref{eq:ratio} holds for $L_f$: $\iota$ multiplies the cubic form of Theorem~\ref{thm:cubic} by $-1$, and $[\sigma_0]\equiv1\pmod{\mathfrak m}$. Complex conjugation, which acts by $(T_c,T_a)\mapsto(T_c,-T_a)$, also preserves \eqref{eq:ratio}, but it does not relate the two interpolation formulas: $L(f,\mathcal W^c,1)=L(f,\mathcal W,1)$, which differs from $L(f,\bar{\mathcal W},1)$ in general.
On the cyclotomic line, Perrin-Riou's factorisation \cite[(1.1)]{PR87} reads
\begin{equation}\label{eq:PR11}
L_f(T_c,0)=L_p(f)\,L_p(f^\chi)\Bigl(\frac{\Omega_f}{\sqrt{d}}\Bigr)^{-1},
\end{equation}
where $L_p(f)$ and $L_p(f^\chi)$ are the Mazur--Swinnerton-Dyer measures with the periods inside. With the normalisation of Lemma~\ref{lem:BK} this gives
\begin{equation}\label{eq:T3}
[T_c^3]L_f=\kappa_F\,c_2(E)\,c_1(E^\chi),\qquad\kappa_F=\frac{\sqrt d\,\Omega^+_E\,\Omega^+_{E^\chi}}{\Omega_f},
\end{equation}
with $c_1(E^\chi)$ the linear coefficient of $L_p(E^\chi,T)$. By Rubin's rank-one theorem \cite[Thm.~8.2]{Rubin92}, \eqref{eq:rubin-rank1}, applied to $E^\chi$ and translated by Lemma~\ref{lem:BK},
\begin{equation}\label{eq:c1}
c_1(E^\chi)=(1-\alpha^{-1})^2\,q_\infty(E^\chi)\,\frac{h_3}{\log_p(1+p)},\qquad q_\infty(E^\chi)=\frac{L'(E^\chi,1)}{\Omega^+_{E^\chi}\,\hat h(P_3)}\in\Q^\times
\end{equation}
($E^\chi$ has the same $a_p$, as $p$ splits in $F$). Here $h_3$ is the height of $P_3$ with the sign of $h^+$ (\S\ref{ss:heights}). In rank one the sign of the height matters. \eqref{eq:c1} is Rubin's theorem \eqref{eq:rubin-rank1}, whose height is Perrin-Riou's and equals $h^+$ on rational points; with $-h_3$ in place of $h_3$ it would be false. We do not use the rank-one case of the formula of \cite{MTT}, whose comparison with \eqref{eq:c1} depends on their normalisations of $L_p$, of the variable and of the pairing; the conventions used here, which are also those of PARI/GP, are confirmed numerically: the rank-one calibration of \S\ref{ss:rankone} checks \eqref{eq:rubin-rank1} with $h^+$ to precision $5^7$; and for $E_2$, of rank one, \texttt{ellpadicregulator} of PARI/GP equals $h^+$ (to $O(5^{12})$ and $O(13^{10})$), and the quotient returned by \texttt{ellpadicbsd}, which conjecturally equals $\#\Sha$ times that regulator, divided by it is $1+O(5^3)$ at $p=5$ and $1+O(13^2)$ at $p=13$. This sign is used in the proof of Theorem~\ref{thm:heegner}. The rationality of $q_\infty(E^\chi)$ follows from the Gross--Zagier formula, and $\kappa_F\deg\varphi\in\Q^\times$ for a modular parametrisation $\varphi\colon X_0(N_E)\to E$, by the relation between $\Omega_f$, $\deg\varphi$ and the area of $E(\C)$ and the relation between the periods of $E$ and of its quadratic twist.

\subsection{The derived Heegner class}
Let $y_s\in E(H_s)\otimes\Zp$ be the regularised Heegner points of conductor $p^s$ of \cite[\S0]{Howard05}, $y_s=\varphi_*(z_s)$ with $z_s=\alpha^{-s}h_s-\alpha^{-s-1}h_{s-1}$ for $s>0$ and $z_0=u_F^{-1}(1-\sigma_\fp/\alpha)(1-\sigma_{\bar\fp}/\alpha)h_0$, which are norm-compatible, and $\kappa_\infty\in H^1_{\rm Iw}(F^a_\infty,T_pE)$ the class of their traces to the anticyclotomic tower. Its image at the bottom is an Euler factor times the Heegner point $y_F$, which is torsion by Gross--Zagier because $L'(E/F,1)=0$; as $E(F)[p]=0$, $\kappa_\infty(0)=0$. Since $E(F^a_\infty)[p]=0$, the map $H^1_{\rm Iw}/T_a\to H^1(F,T_pE)$ is injective and $H^1_{\rm Iw}$ has no $\Lambda$-torsion, so
\[\kappa_\infty=T_a\,\kappa',\qquad \mathcal D:=\kappa'(0)\in H^1_f(F,T_pE)=E(F)\otimes\Zp ,\]
where $\kappa'$ lies in the compact Selmer module by saturation, and the last equality uses the finiteness of $\Sha(E/F)[p^\infty]$. The class $\mathcal D$ is orthogonal to $E(F)$ for the anticyclotomic height. Indeed, $\kappa'$ lifts $\mathcal D$ to the Selmer complex of \S\ref{sec:heegner} modulo $T_a^2$, so the anticyclotomic Bockstein map vanishes on $\mathcal D$, and the anticyclotomic height, which is the pairing defined by this Bockstein map \cite[11.1.4, 11.3.9]{Nekovar06} (up to sign), satisfies $h^F_a(\mathcal D,y)=0$ for all $y\in E(F)$; there is no extended local term at the origin because $H^0(F_v,V^-)=0$ for $v\mid p$, as $\alpha\ne1$. (For the anticyclotomic tower of $K$ and universal norms in both arguments, the corresponding statement is \cite[Appendix, Lemma~A.4]{AH06}.) Writing $\mathcal D=aP_1+bP_2+cP_3$, orthogonality to $P_1,P_2$ gives $c\,v=0$, so $c=0$ by (H6), and orthogonality to $P_3$ gives $(a,b)\parallel(v_2,-v_1)$. Hence
\[\mathcal D=\zeta\,Q_v\qquad\text{for a unique }\zeta\in\Qp .\]

\begin{theorem}[Howard {\cite[Thm.~B]{Howard05}}]\label{thm:howard}
Assume $-d$ odd, $d\ne3$, the Heegner hypothesis and $p$ split in $F$ and ordinary for $f$. Write $L_f=L_{f,0}+L_{f,1}(\gamma_0-1)+\cdots$ along the cyclotomic variable. Then $\kappa\cdot\log_p(\gamma_0)\cdot L_{f,1}=L_{\rm Heeg}$, where $\kappa$ is an Artin symbol and $L_{\rm Heeg}=\varprojlim\sum_\sigma\langle y_s^\vee,y_s^\sigma\rangle_{E,H_s}\,\sigma$, with $y_s=\deg\varphi\cdot y_s^\vee$ and the canonical $p$-adic height of \cite[(9)]{Howard05}.
\end{theorem}

The hypotheses are (H2)--(H3); note that $d\ne3$ because $2$ splits in $F$, so that $-d\equiv1\pmod8$. Theorem~B does not exclude complex multiplication (only the divisibility in \cite[(2)]{Howard05} does). Howard's height is the norm-adapted height, which by \cite[Prop.~3.2.1]{Howard05} and \cite[6.11, 8.1]{Nekovar93} is the Bockstein height of Theorem~\ref{thm:cubic}.

\subsection{The identity}
\begin{lemma}\label{lem:howard-sign}
Let $\langle\ ,\ \rangle^{\rm Howard}$ be the pairing \cite[(9)]{Howard05} on the abelian variety, for Howard's id\`ele class character $\rho$ (normalised by arithmetic reciprocity) and the canonical principal polarisation $\phi_\Theta$. For $P,Q\in E(\Q)\otimes\Qp$, the pairing over $F$ satisfies
\[\langle\phi_\Theta(P),Q\rangle^{\rm Howard}=h^F_c(P,Q)=2\,h^+(P,Q).\]
\end{lemma}
\begin{proof}
Howard's character satisfies $\rho_{\Q,v}(\pi_v)=\log_p\mathrm N(v)$ for $v\nmid p$, the reciprocity map being arithmetic \cite[\S3.3]{Howard05}; as $\rho$ is trivial on $\Q^\times$, $\rho=-\ell$ for the character $\ell$ of Balakrishnan and Besser, with $\ell_p=\log_p$ and $\ell_q(q)=-\log_p q$ \cite[\S4]{BB15}. The abelian pairing is the negative of the pairing on divisors under the principal polarisation \cite[Rem.~3.3.1]{Howard05}, \cite[Prop.~6.9, Thm.~6.10]{BMS26}, and the latter, for the unit-root splitting, is the Coleman--Gross height \cite[Rem.~6.11]{BMS26}, \cite[Thm.~3.6]{IW03}. Hence $\langle\phi_\Theta(P),Q\rangle^{\rm Howard}=-h^{\rm CG}_\rho(P,Q)=h^{\rm CG}_\ell(P,Q)$, and on the diagonal $h^{\rm CG}_\ell(P,P)=-2\log_p(\sigma_p(P)/d_0(P))$ with $x(P)=a/d_0^2$ \cite[Cor.~4.2, 4.3]{BB15}, which is $h^+(P)$ by \S\ref{ss:heights}. Polarisation extends this to pairs, and extending the character to $F$ by the norm multiplies its values on rational points by $[F:\Q]=2$.
\end{proof}

Let
\[S_F=\tfrac12\,q_\infty(E^\chi)\,\kappa_F\deg\varphi\ \in\Q^\times .\]

\begin{theorem}\label{thm:heegner}
Under Hypotheses~\ref{hyp}, $\zeta\ne0$ and
\begin{equation}\label{eq:heegner}
q_p=\frac{1}{S_F}\Bigl(\frac{2\,\ell_a(\sigma_a)\,\zeta}{(1-\alpha^{-1})^{2}}\Bigr)^{2}.
\end{equation}
The right-hand side does not depend on the choice of $\sigma_a$ and $\ell_a$.
\end{theorem}
\begin{proof}
Project Theorem~\ref{thm:howard} to $F^a_\infty$. Since $\kappa_\infty=T_a\kappa'$ and the $\Lambda$-adic height is sesquilinear, $h_\infty(T_a\kappa',T_a\kappa')=T_aT_a^\iota h_\infty(\kappa',\kappa')$ with $T_a^\iota=-T_a+O(T_a^2)$; as $L_{\rm Heeg}$ lies in $(T_a^2)$, the Artin symbol $\kappa$ only contributes through its augmentation $1$, and
\[\begin{split}[T_cT_a^2]\,L_f&=-\frac{h^F_c(\mathcal D,\mathcal D)}{\deg\varphi\,\log_p(\gamma_0)}=-\frac{\zeta^2\,h^F_c(Q_v,Q_v)}{\deg\varphi\,\log_p(1+p)}\\&=\frac{\zeta^2\,Q_2}{\deg\varphi\,\log_p(1+p)},\end{split}\]
where the sign comes from $T_aT_a^\iota=-T_a^2+O(T_a^3)$, and Howard's pairing is $h^F_c$ by Lemma~\ref{lem:howard-sign}. On the other hand, by \eqref{eq:ratio} for $L_f$, \eqref{eq:T3}, \eqref{eq:qp} and \eqref{eq:c1},
\[\begin{split}[T_cT_a^2]\,L_f&=\frac{\log_p(1+p)^2\,Q_2}{\ell_a(\sigma_a)^2\,8h_3\det H}\;\kappa_F\;q_p\frac{(1-\alpha^{-1})^2\det H}{\log_p(1+p)^2}\\&\qquad\times(1-\alpha^{-1})^2q_\infty(E^\chi)\frac{h_3}{\log_p(1+p)}\\&=\frac{\kappa_F\,q_\infty(E^\chi)\,(1-\alpha^{-1})^4\,q_p\,Q_2}{8\,\ell_a(\sigma_a)^2\log_p(1+p)}.\end{split}\]
Equating, and dividing by $Q_2\ne0$, gives $\zeta^2=\kappa_F\deg\varphi\,q_\infty(E^\chi)(1-\alpha^{-1})^4q_p/(8\ell_a(\sigma_a)^2)$, which is \eqref{eq:heegner}. As $q_p\ne0$ by (H5), $\zeta\ne0$. Replacing $\sigma_a$ by $\sigma_a^{w}$, $w\in\Zp^\times$, divides $\mathcal D$ by $w$ and multiplies $\ell_a(\sigma_a)$ by $w$; replacing $\ell_a$ by a multiple $c\ell_a$ multiplies $v$, hence $Q_v$, by $c$ and divides $\zeta$ by $c$.
\end{proof}

\begin{remark}\label{rem:sign}
Two signs enter the comparison of Howard's height with ours, and they cancel: the arithmetic normalisation of the reciprocity map, which makes Howard's character the negative of the one of Coleman and Gross, and the principal polarisation, which makes the abelian pairing the negative of the pairing on divisors. The third sign, $T_aT_a^\iota=-T_a^2$, comes from sesquilinearity, and the global sign of the Bockstein comparison cancels in \eqref{eq:ratio}, as does the inversion $\iota$ that may relate $L_f$ to $L_{\mathcal K}$. The identity therefore carries no undetermined sign. The remaining normalisation is that of \eqref{eq:T3}, which rests on Perrin-Riou's factorisation \cite[(1.1)]{PR87} and Lemma~\ref{lem:BK}. The two examples of \S\ref{ss:heegner-num}, where $q_5$ is known independently, confirm \eqref{eq:heegner}, sign included, modulo $25$ and modulo $5$.
\end{remark}

\section{\texorpdfstring{The square class of $q_p$ and the Birch and Swinnerton-Dyer quotient over $F$}{The square class of q\_p and the BSD quotient over F}}\label{sec:squares}

\subsection{The Manin constant and the normalised quantities}
Let $c_E$ be the Manin constant of the modular parametrisation $\varphi$, so that $\varphi^*\omega_E=c_E\cdot2\pi i f(z)\,dz$ for a N\'eron differential $\omega_E$ (this is the normalisation of \cite[\S2]{ARS06}). Integrating against the complex conjugate gives $c_E^2\,4\pi^2\|f\|^2=\deg\varphi\,\mathrm{vol}(E(\C))$, hence
\[\Omega_f=8\pi^2\|f\|^2=\frac{2\deg\varphi\,\mathrm{vol}(E(\C))}{c_E^2},\qquad c_E^{-2}\,\kappa_F\deg\varphi=\frac{\sqrt d\,\Omega^+_E\Omega^+_{E^\chi}}{2\,\mathrm{vol}(E(\C))}.\]
We put
\[s_F=c_E^{-2}S_F,\qquad \mathcal D^0=c_E^{-1}\mathcal D=\zeta^0Q_v,\qquad \zeta^0=c_E^{-1}\zeta;\]
$\mathcal D^0$ is the derived class built from the Heegner points normalised by $2\pi i f(z)\,dz$, that is, by the elliptic logarithm $\sum a_nn^{-1}q^n$. In these terms \eqref{eq:heegner} reads
\begin{equation}\label{eq:heegner0}
q_p=\frac{1}{s_F}\Bigl(\frac{2\,\ell_a(\sigma_a)\,\zeta^0}{(1-\alpha^{-1})^{2}}\Bigr)^{2},
\end{equation}
and the Manin constant disappears from everything that follows.

\begin{lemma}\label{lem:kappa}
Under (H3), $c_E^{-2}\kappa_F\deg\varphi=\sqrt d\,\Omega^+_E\Omega^+_{E^\chi}/(2\,\mathrm{vol}(E(\C)))$ equals $2$ if $D>0$ and $1$ if $D<0$. Hence $s_F=q_\infty(E^\chi)$ if $D>0$ and $s_F=\tfrac12q_\infty(E^\chi)$ if $D<0$; in particular $s_F/q_\infty(E^\chi)$ is a $p$-adic unit.
\end{lemma}
\begin{proof}
The primes dividing $d$ ramify in $F$ and those dividing $N_E$ split, so $\gcd(d,D)=1$ and $E^\chi=E_{d^2D}$ is given by the global minimal model $y^2=x^3-d^2Dx$. The isomorphism $(X,Y)\mapsto(-dX,d\sqrt{-d}\,Y)$ from $E$ to $E^\chi$ over $F$ shows that the lattice of $dx/2y$ on $E^\chi$ is $\Lambda_E/\sqrt{-d}$, with $\Lambda_E=\Omega_\infty\OK$ as in \S\ref{sec:setup}. If $D>0$, $\Omega_\infty=\omega_1$ is real, both curves have two real components, $\Omega^+_E=2\omega_1$, $\Omega^+_{E^\chi}=2\omega_1/\sqrt d$ and $\mathrm{vol}(E(\C))=\omega_1^2$, so the ratio is $2$. If $D<0$, $\Omega_\infty=\tfrac{1+i}2\omega_1$, both curves have one real component, the least positive real element of $\Omega_\infty\OK$ is $\omega_1$, so $\Omega^+_E=\omega_1$, $\Omega^+_{E^\chi}=\omega_1/\sqrt d$ and $\mathrm{vol}(E(\C))=\omega_1^2/2$, and the ratio is $1$.
\end{proof}

\subsection{The square class}
\begin{theorem}\label{thm:squares}
Under Hypotheses~\ref{hyp},
\[q_p\,s_F\in(\Qp^\times)^2,\qquad s_F=q_\infty(E^\chi)\ \text{or}\ \tfrac12q_\infty(E^\chi)\in\Q^\times .\]
Equivalently $q_pS_F\in(\Qp^\times)^2$.
\end{theorem}
\begin{proof}
In \eqref{eq:heegner0}, $\zeta^0$, $\ell_a(\sigma_a)$ and $1-\alpha^{-1}$ are non-zero elements of $\Qp$, and $S_F/s_F=c_E^2$ is a square.
\end{proof}

Theorem~\ref{thm:squares} is unconditional given the hypotheses; it says nothing about $q_\infty$. To compare with the complex side we need the Birch and Swinnerton-Dyer quotient of $E$ over $F$.

\subsection{\texorpdfstring{The Birch and Swinnerton-Dyer quotient over $F$}{The Birch and Swinnerton-Dyer quotient over F}}
\begin{theorem}\label{thm:bsdF}
Assume (H1) and (H3). Let $2^{j}=[E(F):E(\Q)+E^\chi(\Q)+E(F)_\tors]$ and define $\Shan(E/F)\in\R$ by
\begin{equation}\label{eq:bsdF}
q_\infty(E)\,s_F=4^{1-j}\,\frac{\Shan(E/F)\,\prod_vc_v(E/F)}{\#E(F)_\tors^2}.
\end{equation}
Then $\Shan(E/F)$ is the analytic order of $\Sha(E/F)$, i.e.\ the value predicted for $\#\Sha(E/F)$ by the Birch and Swinnerton-Dyer formula for $E/F$, and $\prod_vc_v(E/F)=\prod_{q\mid N_E}c_q(E)^2$ is a square.
\end{theorem}
\begin{proof}
$L(E/F,s)=L(E,s)L(E^\chi,s)$ vanishes to order three with leading coefficient $\tfrac12L''(E,1)L'(E^\chi,1)$. The period of $E/F$ is $2\,\mathrm{vol}(E(\C))$, the discriminant term is $\sqrt d$, and the regulator of $E(F)$ for the N\'eron--Tate height over $F$ (twice that over $\Q$) is $8\Reg_\infty\hat h(P_3)/4^{j}$, because $E(\Q)$ and $E^\chi(\Q)$ are orthogonal and of index $2^j$ in $E(F)$ modulo torsion. Substituting $\tfrac12L''(E,1)=q_\infty\Omega_E\Reg_\infty$, $L'(E^\chi,1)=q_\infty(E^\chi)\Omega^+_{E^\chi}\hat h(P_3)$ and $s_F=\tfrac12q_\infty(E^\chi)\sqrt d\,\Omega^+_E\Omega^+_{E^\chi}/(2\,\mathrm{vol}(E(\C)))$, the Birch and Swinnerton-Dyer formula for $E/F$ becomes \eqref{eq:bsdF} with $\#\Sha(E/F)$ in place of $\Shan(E/F)$. Every bad prime splits in $F$ by (H3), into two places with the Tamagawa number of $E/\Q_q$.
\end{proof}

The Manin constant does not enter \eqref{eq:bsdF}; with $S_F$ in place of $s_F$ the right-hand side acquires the factor $c_E^2$. Theorem~\ref{thm:bsdF} is an identity of real numbers. Since $q_\infty(E)$ is not known to be rational, it gives no information in $\Qp$ by itself.

\subsection{The rational square class}
The rational number $s_F$ is determined modulo squares by the Tamagawa numbers of $E$. The proof compares the auxiliary twist of rank one with a real twist of rank zero through a genus character, and uses the Birch and Swinnerton-Dyer formula in analytic rank zero for CM curves.

\begin{theorem}\label{thm:tamsquare}
Assume (H1) and (H3). Then
\[\Tam(E)\,s_F\in\Q^{\times2},\qquad\text{equivalently}\qquad \Tam(E)\,S_F\in\Q^{\times2}.\]
More precisely, let $D_0>1$ be a fundamental discriminant such that every prime dividing $2N_Ed$ splits in $\Q(\sqrt{D_0})$ and $L(E^{(D_0)},1)\ne0$, where $E^{(D_0)}$ is the twist of $E$ by $\Q(\sqrt{D_0})$; such $D_0$ exist. Put $q_0=L(E^{(D_0)},1)/\Omega^+_{E^{(D_0)}}$. Then there is $a_F\in\Q^\times$ with
\begin{gather*}
s_F\,q_0=(a_F/c_E)^2,\qquad \Tam(E^{(D_0)})\in\Tam(E)\cdot\Q^{\times2},\\
q_0=\frac{\#\Sha(E^{(D_0)}/\Q)\,\Tam(E^{(D_0)})}{\#E^{(D_0)}(\Q)_\tors^2}.
\end{gather*}
\end{theorem}

\begin{proof}
\emph{Existence of $D_0$.} The sign of the functional equation of $f$ is $+1$, as $\ord_{s=1}L(E,s)=2$. By \cite[Thm.~(ii)]{BFH90}, applied to the set of primes dividing $2N_Ed$, there is a fundamental discriminant $D_0>0$ such that all these primes split in $\Q(\sqrt{D_0})$ and $L(1,f,\chi_{D_0})\ne0$; and $D_0\ne1$ since $L(E,1)=0$. In particular $D_0\equiv1\pmod8$ and $\gcd(D_0,2N_Ed)=1$.

\emph{The genus point.} Let $M=\Q(\sqrt{D_0D_F})$, $D_F=-d$, and let $\eta$ be the genus character of $\mathrm{Pic}(\mathcal O_M)$ attached to the factorisation $D_0D_F$, so that $\mathrm{Ind}_M^\Q\eta=\chi_{D_0}\oplus\chi_F$ and $L(E/M,\eta,s)=L(E^{(D_0)},s)L(E^\chi,s)$. Every prime dividing $N_E$ splits in $\Q(\sqrt{D_0})$ and in $F$, hence in $M$, and $\gcd(D_0d,N_E)=1$. With $\varphi(\infty)=O$, the explicit Gross--Zagier formula of Cai, Shu and Tian \cite[Thm.~1.1]{CST14} (conductor $c=1$, $u=1$ and $\mu(N,D)=0$ here) gives, for $P_\eta=\sum_{[\mathfrak a]}\eta([\mathfrak a])\varphi(P_{\mathfrak a})$,
\[L(E^{(D_0)},1)\,L'(E^\chi,1)=\frac{\Omega_f}{\sqrt{D_0d}\,\deg\varphi}\,\hat h_M(P_\eta).\]
The point $P_\eta$ lies in the $\eta$-part of $E(\Q(\sqrt{D_0},\sqrt{D_F}))\otimes\Q$, which is $(E^{(D_0)}(\Q)\oplus E^\chi(\Q))\otimes\Q$; since $E^{(D_0)}(\Q)$ is finite, $P_\eta=a_FP_3$ with $a_F\in\Q$, and $a_F\ne0$ because the left-hand side is non-zero. The height over $M$ is twice the height over $\Q$ used in the Birch and Swinnerton-Dyer formula \cite[(6.4)]{GZ86}, so $\hat h_M(P_\eta)=2a_F^2\hat h(P_3)$.

\emph{The identity $s_Fq_0=(a_F/c_E)^2$.} By the definitions,
\[s_F=\frac{\sqrt d\,\Omega^+_E\deg\varphi\,L'(E^\chi,1)}{2\,c_E^2\,\Omega_f\,\hat h(P_3)},\qquad\text{so}\qquad s_Fq_0=\Bigl(\frac{a_F}{c_E}\Bigr)^2\frac{\Omega^+_E}{\sqrt{D_0}\,\Omega^+_{E^{(D_0)}}}\]
by the formula above. As $\gcd(D_0,2D)=1$ and $D_0d$ is coprime to the fourth-power-free $D$, the model $y^2=x^3-DD_0^2x$ of $E^{(D_0)}$ is minimal, and $(X,Y)\mapsto(D_0X,D_0^{3/2}Y)$ is an isomorphism over $\R$ carrying its differential $dx/2y$ to $D_0^{-1/2}\,dX/2Y$; so $\Omega^+_{E^{(D_0)}}=\Omega^+_E/\sqrt{D_0}$, with the same number of real components.

\emph{The value of $q_0$.} The twist $E^{(D_0)}$ has complex multiplication by $\Z[i]$ and $L(E^{(D_0)},1)\ne0$, so the Birch and Swinnerton-Dyer formula holds for it, including the prime $2$, by Burungale and Flach \cite[Cor.~2 and the remark following it]{BF24}; $\#\Sha(E^{(D_0)}/\Q)$ is a square by the Cassels--Tate pairing.

\emph{Tamagawa numbers.} At a prime $q\mid N_E$, $D_0$ is a square in $\Q_q$ ($q$ odd splits in $\Q(\sqrt{D_0})$, and $D_0\equiv1\pmod8$), so $E^{(D_0)}\simeq E$ over $\Q_q$ and $c_q(E^{(D_0)})=c_q(E)$. At a prime $\ell\mid D_0$, $E$ has good reduction, $E^{(D_0)}$ has reduction of type $I_0^*$, and by Tate's algorithm \cite[Thm.~4.1, Table~4]{BRSTTW} $c_\ell(E^{(D_0)})=1+\#\{x\in\F_\ell: x^3-Du^2x=0\}$, where $D_0=\ell u$; this is $4$ if $(D/\ell)=1$ and $2$ if $(D/\ell)=-1$. Since $D_0\equiv1\pmod 8$ and $(D_0/q)=1$ for every odd prime $q\mid D$, quadratic reciprocity gives $\prod_{\ell\mid D_0}(D/\ell)=1$, so the number of factors equal to $2$ is even and $\prod_{\ell\mid D_0}c_\ell(E^{(D_0)})$ is a square. Hence $\Tam(E^{(D_0)})\in\Tam(E)\Q^{\times2}$, so $q_0\in\Tam(E)\Q^{\times2}$ and $s_F\in\Tam(E)\Q^{\times2}$. Finally $S_F=c_E^2s_F$.
\end{proof}

All the constants in this proof matter for the square class. As a check, for each of $E_{56},E_{17},E_{-68},E_{-89},E_{-289},E_{-73}$ we computed $q_0$ for four admissible $D_0$ and found $s_F\,q_0$ to be a square in all $24$ cases (it is $64$, $256$ or $1024$); and in all $79$ pairs of \S\ref{ss:bsdF-num} one has $\Tam(E)\,s_F\in\{4,16,64\}$.

\begin{corollary}\label{cor:squares}
Assume Hypotheses~\ref{hyp}. Then $q_p\in\Tam(E)\cdot(\Qp^\times)^2$. If moreover $\Sha(E/\Q)$ is finite, then $q_p/A_E$ is a square in $\Zp^\times$; in particular its reduction modulo $p$ is a quadratic residue.
\end{corollary}
\begin{proof}
By Theorems~\ref{thm:squares} and~\ref{thm:tamsquare}, $q_ps_F$ and $\Tam(E)s_F$ are squares in $\Qp^\times$. If $\Sha(E/\Q)$ is finite, $\#\Sha(E/\Q)$ is a square by the Cassels--Tate pairing, so $A_E\in\Tam(E)\Q^{\times2}$, and $q_p/A_E\in\Zp^\times$ by \S\ref{sec:katz}.
\end{proof}

\begin{remark}\label{rem:existsF}
The conclusion of Corollary~\ref{cor:squares} does not involve $F$, but its proof needs one field $F$ satisfying Hypotheses~\ref{hyp} for the given $(E,p)$. Conditions (H2) and (H5) concern only $(E,p)$. Fields $F$ satisfying (H1) and (H3) exist: since the sign of the functional equation of $f$ is $+1$, \cite[Thm.~(i)]{BFH90}, applied to the set of primes dividing $2pN_E$, gives an imaginary quadratic field in which all of them split and $L(E^\chi,s)$ has a simple zero at $s=1$. Condition (H4) asks $p\nmid h^-_{\mathcal K}$ and (H6) asks $Q_2\ne0$; we do not know how to guarantee them in general, and we check them in each example.
\end{remark}

This is information on the unit $q_p/A_E$ beyond its valuation, which is all that the main conjecture provides; it should be compared with the results of Dokchitser and Dokchitser on Birch and Swinnerton-Dyer quotients modulo squares \cite{DD10}. The class does not depend on $F$, on $p$, or on the finiteness of $\Sha(E/\Q)$: the square class of $q_p$ is that of $\Tam(E)$, which is the class predicted by Conjecture~\ref{conj:main} together with the Birch and Swinnerton-Dyer conjecture and the Cassels--Tate pairing.
If $q_\infty(E)\in\Q$, Theorem~\ref{thm:bsdF} gives $q_\infty s_F\in\Shan(E/F)\cdot\Q^{\times2}$, so that $q_p/q_\infty\in\Shan(E/F)^{-1}(\Qp^\times)^2$, and $q_\infty\in\Tam(E)\Shan(E/F)\Q^{\times2}$ by Theorem~\ref{thm:tamsquare}. Thus, at the fixed prime $p$, Conjecture~\ref{conj:main} is compatible with these results modulo squares exactly when $\Shan(E/F)$ is a square in $\Qp$; the Birch and Swinnerton-Dyer conjecture and the Cassels--Tate pairing predict the stronger assertion that it is a rational square. This compatibility is not a theorem, since the rationality of $q_\infty$ is part of the conjecture.

\section{The conjecture of Agboola and Castella}\label{sec:AC}

\subsection{\texorpdfstring{The BDP $p$-adic $L$-function and the conjecture}{The BDP p-adic L-function and the conjecture}}
Let $\Gamma^a=\mathrm{Gal}(F^a_\infty/F)$, $\Lambda^a=\Zp[[\Gamma^a]]$, $I$ its augmentation ideal, $\iota$ the involution of $\Lambda^a$ induced by $\sigma\mapsto\sigma^{-1}$, and $W$ as in \S\ref{sec:setup}. Under (H2)--(H3), Bertolini, Darmon and Prasanna \cite[Thm.~5.9, Prop.~5.10]{BDP13} and, as an element of the Iwasawa algebra, Castella--Hsieh \cite{CastellaHsieh18} construct $\Lc_p(f)\in\Lambda^a\hat\otimes W$ whose square
\[L^\BDP_p(f)=\Lc_p(f)^2\]
interpolates the central values of the twists of $f/F$ by a family of anticyclotomic Hecke characters of infinite order \cite[Prop.~3.8]{CastellaHsieh18}. The trivial character lies outside the range of interpolation, and the $p$-adic Gross--Zagier formula \cite[Thm.~5.13]{BDP13}, specialised to weight two in \cite[Thm.~3.12]{BDP12} and normalised as in \cite[(1.1)]{AC21}, gives
\begin{equation}\label{eq:BDP1}
L^\BDP_p(f)(\mathbf 1)=u_F^{-2}c_E^{-2}\Bigl(\frac{1-a_p+p}{p}\Bigr)^2\log_E(z_F)^2,
\end{equation}
with $z_F$ a Heegner point, $c_E$ the Manin constant and $u_F=\#\mathcal O_F^\times/2=1$; in our case $z_F$ is torsion and $L^\BDP_p(f)\in I$.

Agboola and Castella \cite[Conj.~1.1, Conj.~4.2]{AC21} predict the order of vanishing and the leading coefficient of $L^\BDP_p(f)$ at $\mathbf 1$ (we use the numbering of the published version). Write $r^\pm$ for the ranks of the eigenspaces of complex conjugation on $E(F)$; here $(r^+,r^-)=(2,1)$, and since $|r^+-r^-|=1$ their derived regulator equals the regulator below as soon as the conjecture of Mazur and Bertolini--Darmon holds \cite[Conj.~3.8, Lemma~3.10]{AC21}; in our situation it does (Lemma~\ref{lem:MBD} below). Let $P_1,P_2,P_3\in E(F)$ be linearly independent and $t=[E(F):\Z P_1+\Z P_2+\Z P_3]$, an index in $E(F)$ itself, so that torsion contributes to $t$; the anticyclotomic height takes values in $I/I^2$ through $\sigma_a-1\mapsto\ell_a(\sigma_a)$, i.e.\ $h^F_a$ is read as $h^F_a\cdot(\sigma_a-1)/\ell_a(\sigma_a)$. The regulator of \cite[Def.~3.4]{AC21} is
\[\Reg=\frac1{t^2}\sum_{i,k}(-1)^{i+k}P_i\otimes P_k\,R_{ik},\]
$R_{ik}$ the $(i,k)$-minor of the matrix of anticyclotomic heights; each $R_{ik}$ is a sum of products of two heights, an element of $I^2/I^3$.

\emph{The product of two heights.} There are two natural ways to multiply two elements $x,y\in I/I^2$ into $I^2/I^3$: the product $xy$ of the ring $\Lambda^a$, and the sesquilinear product $x\,y^\iota$, which is $-xy$ because $\iota$ acts as $-1$ on $I/I^2$. We write $\Reg^\iota$ for the regulator computed with the sesquilinear product; when $r^+ + r^- = 3$, as here, every term of a minor is a product of two heights, so $\Reg^\iota=-\Reg$. The sesquilinear product is the one that occurs in the $\Lambda$-adic theory: Howard's height is sesquilinear, $h_\infty(T_ax,T_ay)=T_aT_a^\iota\,h_\infty(x,y)$, and this is how the anticyclotomic heights enter the proof of Theorem~\ref{thm:heegner}. In rank one (the case of \cite[Rem.~1.2]{AC21}) no product occurs and the two readings agree. We state the conjecture with $\Reg^\iota$, and justify the choice after Theorem~\ref{thm:AC}.

\begin{conjecture}[Agboola--Castella, analytic form]\label{conj:AC}
$\Shan(E/F)\in\Q$, $L^\BDP_p(f)\in I^2$, and its image $\bar L^\BDP_p$ in $I^2/I^3$ satisfies
\[\bar L^\BDP_p=\Bigl(\frac{1-a_p+p}{p}\Bigr)^2\,(\log_E\otimes\log_E)(\Reg^\iota)\,\Shan(E/F)\,{\textstyle\prod_\ell} c_\ell(E)^2 .\]
\end{conjecture}

Conjecture~\ref{conj:AC} is the square of \cite[Conj.~4.3]{AC21}, with $\#\Sha(E/F)$ replaced by its analytic value $\Shan(E/F)$ of Theorem~\ref{thm:bsdF}, which makes it a statement about $L$-functions and heights only, like Conjecture~\ref{conj:main}. Under the Birch and Swinnerton-Dyer conjecture for $E/F$ the two values of $\Sha$ agree, but the sign does not change: the literal form of \cite[Conj.~4.2]{AC21} differs from Conjecture~\ref{conj:AC} by the sign $\Reg=-\Reg^\iota$ discussed after Theorem~\ref{thm:AC}. Torsion enters through the index $t$, as in the rank-one index formula on which the conjecture is modelled \cite[Rem.~1.2]{AC21} (compare \cite[Conj.~V.2.2 and (2.3), p.~311]{GZ86} and \cite{Jetchev08}). This is invisible up to $p$-adic units, which is the level of the results of \cite{AC21}, but it matters for the exact formula: every curve $E_D$ has the point $(0,0)$ of order two, and the $79$ cases of \S\ref{ss:bsdF-num} show that \eqref{eq:bsdF}, hence the comparison below, requires the factor $\#E(F)_\tors^{-2}$ contained in $t^{-2}$.

\begin{lemma}\label{lem:ACreg}
In the situation of \S\ref{sec:heegner}, with $P_1,P_2$ a basis of $E(\Q)$ and $P_3$ of $E^\chi(\Q)$ modulo torsion, $t=2^j\,\#E(F)_\tors$ and
\[(\log_E\otimes\log_E)(\Reg^\iota)=\frac{\log_E(Q_v)^2}{4^{j}\,\#E(F)_\tors^2}\cdot\frac{(\sigma_a-1)^2}{\ell_a(\sigma_a)^2},\qquad\log_E(Q_v)=v_2\lambda_1-v_1\lambda_2 ,\]
where $\lambda_a=\lambda(P_a)$ are the formal logarithms of \S\ref{sec:rubin}.
\end{lemma}
\begin{proof}
The subgroup $\Z P_1+\Z P_2+\Z P_3$ meets $E(F)_\tors$ trivially and has index $\#E(F)_\tors$ in $E(\Q)+E^\chi(\Q)+E(F)_\tors$, which has index $2^j$ in $E(F)$. In the basis $P_1,P_2,P_3$ the matrix of anticyclotomic heights is $H_a$, of rank two, and $\adj(H_a)=-n\,{}^tn$ with $n=(v_2,-v_1,0)$; hence $\Reg=-t^{-2}\,Q_v\otimes Q_v$ and $\Reg^\iota=t^{-2}\,Q_v\otimes Q_v$, with the heights read in $I/I^2$ through $\ell_a(\sigma_a)^{-1}(\sigma_a-1)$.
\end{proof}

\begin{lemma}\label{lem:MBD}
In the situation of \S\ref{sec:heegner}, let $S_p=E(F)\otimes\Zp$ (Lemma~\ref{lem:consequences}(2)), let $S_p=S^{(1)}_p\supseteq S^{(2)}_p\supseteq\cdots$ be the filtration by the null-spaces of the derived anticyclotomic height pairings $h^{(k)}_p$, $h^{(1)}_p$ being the anticyclotomic height \cite[(3.6)]{AC21}, and let $US_p$ be the submodule of universal norms. Then
\[S^{(k)}_p\otimes\Qp=US_p\otimes\Qp=\Qp\,Q_v\qquad(k\ge2).\]
In particular the conclusion of the conjecture of Mazur and Bertolini--Darmon \cite[Conj.~3.8]{AC21} holds up to finite index, with $e_1=2$ and $e_k=0$ for $k\ge2$, and the derived regulator of \cite[Def.~3.9]{AC21} equals the regulator $\Reg$.
\end{lemma}
\begin{proof}
The pairing $h^{(1)}_p$ is the anticyclotomic $p$-adic height, which agrees with $h^F_a$ up to the identification of $I/I^2$ and a sign \cite[11.3.9]{Nekovar06}, \cite[Rem.~4.3]{CHKLL}; neither affects null-spaces. In the basis $P_1,P_2,P_3$ its matrix is $H_a$, of rank two because $v\ne0$ by (H6), with null-space $\Qp\,Q_v$; so $S^{(2)}_p\otimes\Qp=\Qp\,Q_v$, of rank $|r^+-r^-|=1$. The universal norms lie in the null-space of every $h^{(k)}_p$ \cite[\S3]{AC21}, so $US_p\otimes\Qp\subseteq S^{(k)}_p\otimes\Qp\subseteq\Qp\,Q_v$ for $k\ge2$, and $US_p$ contains $\mathcal D=\zeta Q_v$ with $\zeta\ne0$ (Theorem~\ref{thm:heegner}). Hence all these spaces are $\Qp\,Q_v$. The proof of \cite[Lemma~3.10]{AC21} uses only that universal norms are in the null-space of $h^{(1)}_p$ and that $S^{(2)}_p=S^{(3)}_p=\cdots=US_p$; since the regulators are taken after tensoring with $\Q$, the equalities after $\otimes\Qp$ suffice, and it gives the last assertion.
\end{proof}

The same computation in the basis of \cite[Conj.~1.1]{AC21} (a basis $x_1,x_2$ of the kernel of localisation at $\fp$, for instance $x_1=Q_\lambda$ and $x_2=\lambda_3P_1-\lambda_1P_3$, completed by $y_p$) gives the same value: there $\log_E(y_p)^2\det\bigl(h^F_a(x_i,x_k)\bigr)t_M^{-2}=-\log_E(Q_v)^2/t^2$ for the ring product. The minus sign is intrinsic: $x_1$ is isotropic, so the anticyclotomic form on the kernel of localisation is hyperbolic; compare \cite[Lemma~3.14]{AC21}.

The condition $\log_E(Q_v)\ne0$ says that $Q_v$ is not proportional to $Q_\lambda$, i.e.\ that $h_a^F(Q_\lambda,P_3)\ne0$. It holds in our examples: at $p=5$ the ratios $v_2/v_1$ and $\lambda_2/\lambda_1$ differ modulo $25$ for $E_{56}$ and modulo $5$ for $E_{-289}$.

\subsection{Explicit reciprocity}
We work in a semilocal Iwasawa algebra and project to $\Lambda^a$ only at the end. Let $H_{p^\infty}$ be the union of the ring class fields of $F$ of $p$-power conductor, $\tilde\Gamma=\mathrm{Gal}(H_{p^\infty}/F)$, $\pi\colon\tilde\Gamma\to\Gamma^a$ the projection, whose kernel is finite, and $W'$ the extension of $W$ generated by the values of the character $\psi$ below.

\emph{The auxiliary character.} Let $\psi$ be a Hecke character of $F$ of infinity type $(1,-1)$ and conductor $1$: on principal ideals put $\psi((a))=\bar a/a$, which is well defined because $\mathcal O_F^\times=\{\pm1\}$ and is the value forced by $\psi_\infty(z)=z/\bar z$ and the triviality of $\psi$ on $F^\times$ (with the conventions of \cite[\S3.3]{CastellaHsieh18} and \cite[\S2.6]{Castella20}), and extend it to the ideal group, which is possible since the class group is finite and $\C^\times$ is divisible. It is anticyclotomic, because $\psi(\mathfrak a)\psi(\bar{\mathfrak a})=\psi((\mathrm N\mathfrak a))=1$. Its $p$-adic avatar $\hat\psi$ is a character of $\tilde\Gamma$, crystalline at the places above $p$, which in general does not factor through $\Gamma^a$; so $\psi$ satisfies the hypotheses of \cite[Thm.~5.7]{CastellaHsieh18}, with $c_o=1$.

\emph{The classes.} By \cite[Def.~5.2]{CastellaHsieh18} with $r=1$, the class $z_{f,1,\alpha}$ of \cite[(5.1)]{CastellaHsieh18} is $\varprojlim_n\alpha^{-n}(z_{p^n}-\alpha^{-1}z_{p^{n-1}})$, where in weight two the classes $z_c$ are the Kummer images of the Heegner points of conductor $c$ on $J_0(N_E)$ \cite[Rem.~6.6]{Castella20}. These are the regularised points $z_s$ of \cite[\S0]{Howard05}, so the whole Iwasawa class, and not only its bottom class, which vanishes here, is Howard's: the images of $z_f=z_{f,1,\alpha}$ under $\varphi_*$ and corestriction to the anticyclotomic tower of $F$ give $\kappa_\infty$. (Castella's theorem \cite[Thm.~6.5, Rem.~6.6]{Castella20} identifies the same class with the weight-two specialisation of Howard's big Heegner points \cite{Howard07}.)

\emph{The reciprocity law over $\tilde\Gamma$.} The explicit reciprocity law of Castella and Hsieh \cite[Thm.~5.7]{CastellaHsieh18} for $\psi$ reads, for $r=1$, where its factor $-c_o^{\,r-1}$ equals $-1$,
\[\bigl\langle\mathfrak L_{\fp,\psi}(z_f),\,\omega_f\otimes t^{-1}\bigr\rangle=-\Lc_{\fp,\psi}(f)\cdot\sigma_{-1,\fp}\qquad\text{in }W'[[\tilde\Gamma]],\]
where $\sigma_{-1,\fp}=\mathrm{rec}_\fp(-1)$ has order two, $\mathfrak L_{\fp,\psi}(z_f)$ is the big logarithm of the twist $z_f\otimes\hat\psi^{-1}$, read on $\tilde\Gamma$ \cite[Def.~5.6]{CastellaHsieh18}, and $\Lc_{\fp,\psi}(f)$ is the function of \cite[Prop.~3.8]{CastellaHsieh18}. The auxiliary period does not survive: the vector $\omega_f\otimes t^{-1}$ is defined through $\langle x,\omega_{f,\psi}\rangle=\Omega_\psi\langle x,\omega_f\otimes t^{-1}\rangle$ \cite[(5.3)]{CastellaHsieh18}, and in the proof of \cite[Thm.~5.7]{CastellaHsieh18} the period $\Omega_\psi$ of $\hat\psi|_{G_{F_\fp}}$ cancels against the one in the interpolation formula of $\Lc_{\fp,\psi}(f)$ (equations (5.6)--(5.7) of loc.\ cit.).

\emph{Recentring and projection.} For a character $\eta$ of $\tilde\Gamma$, let $\mathrm{Tw}_\eta$ be the automorphism $g\mapsto\eta(g)g$ of $W'[[\tilde\Gamma]]$, so that $(\mathrm{Tw}_\eta M)(\chi)=M(\eta\chi)$, and let $\iota$ be the involution $g\mapsto g^{-1}$. In \cite{CastellaHsieh18} the value of $\Lc_{\fp,\psi}(f)$ at $\hat\phi^{-1}$ corresponds to the character $\phi$, and the trivial character to $\phi=\psi$ (see (5.6) of loc.\ cit.). So the function centred at the trivial character is
\[\Lc_p(f)=\pi\bigl(\iota\,\mathrm{Tw}_{\hat\psi^{-1}}\Lc_{\fp,\psi}(f)\bigr)\in\Lambda^a\hat\otimes W',\qquad\Lc_p(f)(\chi)=\Lc_{\fp,\psi}(f)\bigl(\hat\psi^{-1}(\chi\circ\pi)^{-1}\bigr);\]
$L^\BDP_p(f)=\Lc_p(f)^2$, and its value at the trivial character is \eqref{eq:BDP1}. Apply $\pi\circ\iota\circ\mathrm{Tw}_{\hat\psi^{-1}}$ to both sides of the reciprocity law. On the right, $\sigma_{-1,\fp}$ becomes $\pm\pi(\sigma_{-1,\fp})=\pm1$, because $\sigma_{-1,\fp}$ has order two, $\hat\psi(\sigma_{-1,\fp})=\pm1$ and $\Gamma^a$ is torsion-free. On the left, the big logarithm $\mathfrak L$ of \cite[Thm.~3.7]{Castella20} (of which $\mathfrak L_{\fp,\psi}$ is the specialisation, as in the proof of \cite[Thm.~6.5]{Castella20}) is linear over the Iwasawa algebra of the local Galois group, and its specialisation at a character is given by the Bloch--Kato logarithm of the correspondingly specialised class. So the transformed left-hand side specialises at $\chi$ to the logarithm of the $\chi$-specialisation of $\pi_*z_f$: taking $\hat\phi=\hat\psi\chi$ in \cite[(5.6)]{CastellaHsieh18}, the argument is $\hat\psi^{-1}\chi^{-1}$ and the class is $z_f^{\chi}$. This logarithm is computed in $D_{\rm dR}(F^+V_f(1)\otimes\hat\psi^{-1}\otimes\hat\psi)=D_{\rm dR}(F^+V_f(1))$ and paired with $\omega_f$; the identification of these two Dieudonn\'e modules is multiplication by the period of $\hat\psi|_{G_{F_\fp}}$, which is the factor normalised away in \cite[(5.3)]{CastellaHsieh18}. Hence, with $\mathfrak L$ formed with $\omega_f$, $\mathfrak L(\pi_*z_f)$ equals $\Lc_p(f)$ up to a sign and the involution $\iota$, and $\iota$ acts as $-1$ on $I/I^2$. The big logarithm on $T_f$ is formed with $\omega_f=2\pi if(z)\,dz$; since $\varphi^*\omega_E=c_E\,\omega_f$ (\S\ref{sec:squares}), transport by $\varphi_*$ to $T_pE$, with the logarithm $\log_E$ of $\omega_E$, divides it by $c_E$. So, writing $\mathfrak L$ also for the big logarithm on $T_pE$ formed with $\omega_E$, $\mathfrak L(c_E^{-1}\kappa_\infty)$ and $\Lc_p(f)$ have the same image in $I/I^2$ up to sign.

\emph{The leading term.} By \cite[Thm.~3.7]{Castella20}, at the trivial character of $\Gamma^a$ (Hodge--Tate weight $w=1$, conductor $1$, Frobenius eigenvalue $\varphi=\alpha$, the same for $F^+V_f(1)\otimes\hat\psi^{-1}$ at $\hat\psi$ as for $F^+V_f(1)$ at $\mathbf 1$) the big logarithm is the Bloch--Kato logarithm times
\[\Bigl(1-\frac{\varphi}{p^{w}}\Bigr)\Bigl(1-\frac{p^{w-1}}{\varphi}\Bigr)^{-1}=\frac{1-\alpha/p}{1-\alpha^{-1}}=\frac{1-\beta^{-1}}{1-\alpha^{-1}}.\]
This agrees with \eqref{eq:BDP1} and \cite[Thm.~A]{Castella20}, since $(1-\alpha^{-1})(1-\beta^{-1})=(1-a_p+p)/p$. As $\kappa_\infty=T_a\kappa'$ and $\mathfrak L$ is $\Lambda^a$-linear, for a sign $\pm$ independent of $\chi$,
\begin{equation}\label{eq:reciprocity}
\Lc_p(f)\equiv\pm\frac{1-\beta^{-1}}{1-\alpha^{-1}}\,\log_E(\mathcal D^0)\,(\sigma_a-1)\pmod{I^2},\qquad\log_E(\mathcal D^0)=\zeta^0\log_E(Q_v).
\end{equation}
Only the square of \eqref{eq:reciprocity} is used below, so the sign plays no role.

\subsection{The equivalence}
\begin{theorem}\label{thm:AC}
Assume Hypotheses~\ref{hyp} and $\log_E(Q_v)\ne0$. Then $L^\BDP_p(f)\in I^2\setminus I^3$, and
\begin{equation}\label{eq:ACus}
\bar L^\BDP_p=\Bigl(\frac{1-a_p+p}{p}\Bigr)^2\,\frac{\log_E(Q_v)^2}{4\,\ell_a(\sigma_a)^2}\,q_p\,s_F\,(\sigma_a-1)^2 .
\end{equation}
Consequently:
\begin{enumerate}
\item Conjecture~\ref{conj:AC} for $(E,F,p)$ holds if and only if $q_\infty\in\Q$ and $q_p=q_\infty$, that is, if and only if Conjecture~\ref{conj:main} holds for $(E,p)$;
\item the square class of $\bar L^\BDP_p/(\log_E(Q_v)\,(\sigma_a-1)/\ell_a(\sigma_a))^2$ in $\Qp^\times/(\Qp^\times)^2$ is trivial, which is the class predicted by Conjecture~\ref{conj:AC} when $\Shan(E/F)$ is a square in $\Qp$.
\end{enumerate}
\end{theorem}
\begin{proof}
By \eqref{eq:reciprocity}, $\bar L^\BDP_p=\bigl(\tfrac{1-\beta^{-1}}{1-\alpha^{-1}}\bigr)^2(\zeta^0)^2\log_E(Q_v)^2(\sigma_a-1)^2$, and by \eqref{eq:heegner0}, $(\zeta^0)^2=q_ps_F(1-\alpha^{-1})^4/(4\ell_a(\sigma_a)^2)$; this gives \eqref{eq:ACus}, which is non-zero. By Lemma~\ref{lem:ACreg} and Theorem~\ref{thm:bsdF}, the right-hand side of Conjecture~\ref{conj:AC} is
\[\begin{split}&\Bigl(\frac{1-a_p+p}{p}\Bigr)^2\frac{\log_E(Q_v)^2}{4^{j}\ell_a(\sigma_a)^2}\,\frac{\Shan(E/F)\prod_\ell c_\ell^2}{\#E(F)_\tors^2}(\sigma_a-1)^2\\&\qquad=\Bigl(\frac{1-a_p+p}{p}\Bigr)^2\frac{\log_E(Q_v)^2}{4\,\ell_a(\sigma_a)^2}\,q_\infty s_F\,(\sigma_a-1)^2,\end{split}\]
so the two agree exactly when $q_p=q_\infty$, the rationality of $q_\infty$ being equivalent to that of $\Shan(E/F)$. For (2), use Theorem~\ref{thm:squares} and $-1\in(\Qp^\times)^2$.
\end{proof}

\emph{The reading of the products.} With $\Xi=2\ell_a(\sigma_a)\zeta^0/(1-\alpha^{-1})^2$ as in Table~\ref{tab:heegner}, \eqref{eq:reciprocity} and $(1-\alpha^{-1})(1-\beta^{-1})=(1-a_p+p)/p$ give, without using Theorem~\ref{thm:heegner},
\[\bar L^\BDP_p=\Bigl(\frac{1-a_p+p}{p}\Bigr)^2\frac{\log_E(Q_v)^2}{4\,\ell_a(\sigma_a)^2}\,\Xi^2\,(\sigma_a-1)^2 .\]
By Lemma~\ref{lem:ACreg} and Theorem~\ref{thm:bsdF}, Conjecture~\ref{conj:AC} is therefore equivalent to $\Xi^2=q_\infty s_F$, while the same conjecture with the ring product ($\Reg=-\Reg^\iota$) is equivalent to $\Xi^2=-q_\infty s_F$. The two readings can therefore be tested on Heegner points without knowing $q_p$; in \S\ref{ss:readings} we do this for the two pairs of \S\ref{ss:heegner-num}, and both agree with the sesquilinear reading and not with the ring product. This comparison involves only the explicit reciprocity law, the Heegner points and the complex side, not Theorem~\ref{thm:heegner} or $q_p$.

The sign is also visible inside \cite{AC21}. There $L_p(f)=\Lc_p(f)^2$, without involution \cite[\S2]{AC21}, whereas the conjecture of Bertolini and Darmon that \cite[\S3]{AC21} refines concerns $L_K=\theta\theta^\iota$, $\theta$ the Heegner distribution, and $\theta\theta^\iota\equiv(-1)^\rho\theta^2$ in $I^{2\rho}/I^{2\rho+1}$ for $\theta\in I^\rho$. By \cite[Lemma~3.14]{AC21} the first partial regulator is $R^{(1)}=-\det\bigl(h_p(y_i^+,y_j^-)\bigr)^2$, so the square of the conjecture for the square root $\Lc_p(f)$ \cite[Conj.~4.3]{AC21} is the conjecture for $L_p(f)$ \cite[Conj.~4.2]{AC21} with $\Reg_{p,\mathrm{der}}$ replaced by $(-1)^\rho\Reg_{p,\mathrm{der}}$. Here $\rho=\max\{r^+,r^-\}-1=1$, and our Conjecture~\ref{conj:AC} is the square of \cite[Conj.~4.3]{AC21}; the literal form of \cite[Conj.~4.2]{AC21} differs from it by this sign. The formulation of Castella, Hsu, Kundu, Lee and Liu \cite[Conj.~1.1, Def.~3.14]{CHKLL} points the same way. Their partial regulators pair the $\fp$-strict and the $\fpb$-strict Selmer groups, $R^{(i)}=\det h^{(i)}(x_{i,j},y_{i,j'})$, which complex conjugation identifies, and their normalisation $\det(M)^{-2}\Reg_{\fp,\mathrm{der}}$ \cite[Rem.~1.3]{CHKLL} is independent of $M$ when the second basis is the image $y=\tau x$ of the first. In our situation, with $x_1=Q_\lambda$, $x_2=\lambda_3P_1-\lambda_1P_3$ and $s_\fp=P_1$, this gives $\log_\fp(s_\fp)^2\det(M)^{-2}R^{(1)}=+\log_E(Q_v)^2$, the value of $(\log_E\otimes\log_E)(\Reg^\iota)$ in Lemma~\ref{lem:ACreg} up to $t^{-2}$ and the identification of $I/I^2$, and not its negative. Sano's derived Bockstein regulator has the same structure: the Bockstein map induces a pairing between the $\fp$-strict and the $\fpb$-strict Selmer groups, through a duality twisted by $\iota$ \cite[Prop.~4.4, Rem.~4.7]{Sano23}. Both works hold up to $p$-adic units, like \cite{AC21}, so they do not decide the sign. We read this as a sign to be fixed in the transfer from $\theta\theta^\iota$ to $\Lc_p(f)^2$, invisible up to $p$-adic units and in rank one, rather than as evidence against the conjecture of Agboola and Castella. It does not affect the results of \cite{AC21}, which hold up to $p$-adic units, nor their rank-one consistency check \cite[Rem.~1.2]{AC21}, which involves no product of heights. The sign in Theorem~\ref{thm:heegner} itself is determined in Lemma~\ref{lem:howard-sign} and Remark~\ref{rem:sign}.

\subsection{\texorpdfstring{The conjecture up to a $p$-adic unit}{The conjecture up to a p-adic unit}}
Because Theorem~\ref{thm:AC} determines $\bar L^\BDP_p$ exactly, the form of the conjecture of Agboola and Castella that is accessible to Iwasawa theory, the leading-term formula up to a $p$-adic unit with the algebraic $\#\Sha(E/F)$, follows in our situation without the BDP main conjecture: besides Theorem~\ref{thm:AC}, the proof uses only the main conjectures of Rubin and Hsieh, through Sections~\ref{sec:katz} and~\ref{sec:heegner}.

\begin{corollary}\label{cor:ACunit}
Under the hypotheses of Theorem~\ref{thm:AC}, $\bar L^\BDP_p$ equals, up to a unit of $\Zp^\times$, the right-hand side of \cite[Conj.~4.2]{AC21}, that is, of Conjecture~\ref{conj:AC} with $\#\Sha(E/F)$ in place of $\Shan(E/F)$ and with either reading of the products of heights. If $q_\infty\in\Q$ and $\Sha(E/F)$ is finite, the unit is $q_p\,\Shan(E/F)/(q_\infty\#\Sha(E/F))$ for the sesquilinear reading; in general only $\#\Sha(E/F)[p^\infty]$, which is finite, enters.
\end{corollary}
\begin{proof}
By \eqref{eq:ACus} and Lemma~\ref{lem:ACreg}, and since $t=2^j\#E(F)_\tors$ and $\prod_vc_v(E/F)=\prod_qc_q(E)^2$ are prime to $p$ (the Tamagawa numbers of $E_D$ are at most $4$ and the torsion is $2$-primary), the statement amounts to $v_p(q_ps_F)=v_p(\#\Sha(E/F)[p^\infty])$. As $p$ is odd, $\Sha(E/F)[p^\infty]\simeq\Sha(E/\Q)[p^\infty]\oplus\Sha(E^\chi/\Q)[p^\infty]$, and all three groups are finite by Lemma~\ref{lem:consequences}(2). By \S\ref{sec:katz}, $v_p(q_p)=v_p(\#\Sha(E/\Q)[p^\infty])$. By Lemma~\ref{lem:kappa}, $v_p(s_F)=v_p(q_\infty(E^\chi))$. Finally $v_p(q_\infty(E^\chi))=v_p(\#\Sha(E^\chi/\Q)[p^\infty])$: Rubin's theorem \eqref{eq:rubin-rank1} gives $q_p(E^\chi)=q_\infty(E^\chi)$ for the rank-one analogue $q_p(E^\chi)=c_1(E^\chi)\log_p(1+p)/((1-\alpha^{-1})^2h_3)$ of \eqref{eq:qp}, and the argument of \S\ref{sec:katz} in rank one (the main conjecture \cite[Thm.~12.3]{Rubin91} and \cite[Thm.~V.8(3)]{PR84}, with $h_3\ne0$ by Lemma~\ref{lem:consequences}(3)) gives $v_p(q_p(E^\chi))=v_p(\#\Sha(E^\chi/\Q)[p^\infty])$. The two readings of the products differ by $-1$.
\end{proof}

So the conjecture of Agboola and Castella, in the form in which it is stated and studied in \cite{AC21,CHKLL,Sano23}, holds for these CM curves; what Theorem~\ref{thm:AC} adds is the unit itself. When $q_\infty\in\Q$ and $\Sha(E/F)$ is finite, it is $q_p/q_\infty$ times the rational square-class factor $\Shan(E/F)/\#\Sha(E/F)$, its square class is trivial whenever $\Shan(E/F)$ is a square, and the exact formula, in the analytic form of Conjecture~\ref{conj:AC} with the sign corrected, is equivalent to $q_p=q_\infty$.

\subsection{Relation with other results}
Theorem~1.3 of \cite{AC21} proves the leading coefficient formula up to a $p$-adic unit for the characteristic power series of the relevant Selmer group, assuming that $\Sha(E/F)[p^\infty]$ is finite, that $\bar\rho_{E,p}\colon G_\Q\to\mathrm{Aut}(E[p])$ is surjective, that $E[p]$ is ramified at every prime dividing $N_E$ and that $p\nmid\#E(\F_p)$; surjectivity excludes complex multiplication, so for CM curves this result is not available, rather than vacuous. Castella, Hsu, Kundu, Lee and Liu \cite[Thm.~5.1]{CHKLL} prove the algebraic leading-coefficient formula up to a unit under hypotheses that allow complex multiplication; their analytic consequence \cite[Cor.~5.2]{CHKLL} again requires a surjective residual representation. Sano \cite[Thm.~4.13]{Sano23} deduces \cite[Conj.~4.2]{AC21}, up to a unit of $W$, from the anticyclotomic Iwasawa--Greenberg main conjecture for $(E,F,p)$, assuming only the Heegner hypothesis, ordinarity, $E(F)[p]=0$ and the finiteness of $\Sha(E/F)[p^\infty]$ \cite[Hyp.~3.1]{Sano23}; these hold here by Lemma~\ref{lem:consequences}. In our situation the main conjecture is known after inverting $p$ \cite[Thm.~1.2.4(a)]{BCS25}, which with Sano's theorem gives \cite[Conj.~4.2]{AC21} only up to a power of $p$. Corollary~\ref{cor:ACunit} gives it up to a unit directly, without any main conjecture for $L^\BDP_p$, and in turn yields the integral form of the main conjecture (Theorem~\ref{thm:bdpMC}). Note that the main conjecture is not a consequence of Theorem~\ref{thm:planeMC}: its Selmer group, relaxed at one prime of $F$ above $p$ and strict at the other, corresponds over $\mathcal K=KF$ to the $p$-ordinary CM type induced from $F$, whereas Theorem~\ref{thm:planeMC} concerns the ordinary Selmer group, which corresponds to the CM type induced from $K$, and \cite[Thm.~8.16]{Hsieh14} does not apply to the former, since the branch character is ramified at a place of its complement. The proof of the main conjecture in \cite{BCK21}, as used in \cite[Cor.~3.11]{Sano23} and \cite[Cor.~5.2]{CHKLL}, requires a surjective residual representation, and the anticyclotomic main conjecture of \cite{AH06} concerns the CM field $K$ rather than $F$; the input we use is the base change method of Burungale, Castella and Skinner.

\subsection{The anticyclotomic main conjecture}
Let $X_{\rm Gr}$ be the Pontryagin dual of the Greenberg Selmer group of $E[p^\infty]$ over $F^a_\infty$ that is relaxed at the prime $v$ of $F$ above $p$ used to define $L^\BDP_p$ and strict at $\bar v$, and $\Lambda^{a,\rm ur}=\Lambda^a\hat\otimes W$.

\begin{theorem}\label{thm:bdpMC}
Assume (H2) and (H3).
\begin{enumerate}
\item $X_{\rm Gr}$ is $\Lambda^a$-torsion, and $\mathrm{char}(X_{\rm Gr})=(L^\BDP_p(f))$ in $\Lambda^{a,\rm ur}\otimes\Qp$.
\item Under the hypotheses of Theorem~\ref{thm:AC}, the equality holds in $\Lambda^{a,\rm ur}$.
\end{enumerate}
\end{theorem}
\begin{proof}
(1) is \cite[Thm.~1.2.4(a)]{BCS25}. Its hypotheses are: $p>3$ of good ordinary reduction, $E[p]$ irreducible as a $G_\Q$-module, $-d$ odd and $\ne-3$, the Heegner hypothesis and $p$ split in $F$; they hold by (H2), (H3) and Lemma~\ref{lem:consequences}(4), whose proof uses only (H2) and (H3). Complex multiplication is not excluded there: the surjectivity of $\bar\rho_{E,p}$ is required only in part~(b) of that theorem. Since $\bar\rho$ is dihedral, we check the irreducibility conditions used in the proof \cite[Thm.~4.2.1, Prop.~4.2.2, Prop.~5.2.1, Lemma~5.2.3]{BCS25}. Irreducibility over $G_F$ is Lemma~\ref{lem:consequences}(4). Let $F'$ be the auxiliary real quadratic field of \cite[Lemma~5.2.3]{BCS25}: $p$ is inert in $F'$, and $D_{F'}$ is odd and prime to $N_E$, by conditions (i)--(iii) of \cite[Prop.~5.2.1]{BCS25}. For a number field $L\not\ni i$, $\bar\rho|_{G_L}$ is induced from $\bar\chi_\fp|_{G_{L(i)}}$, so it is irreducible if $\bar\chi_\fp\ne\bar\chi_\fpb$ on $G_{L(i)}$.
\begin{itemize}
\item Condition (iv), irreducibility over $M=F'F$: $i\notin M$, since $M/\Q$ is unramified at $2$ ($D_{F'}$ and $d$ being odd), and $p$ is unramified in $M$, so the inertia argument of Lemma~\ref{lem:consequences}(4) applies.
\item Condition (v), irreducibility over $F'(\mu_p)$: $i\notin F'(\mu_p)$, since $2\nmid pD_{F'}$. Inertia at $\fp$ and at $\fpb$ shows that $G_K$ acts on $E[p]=E[\fp]\oplus E[\fpb]$ through all of $\F_p^\times\times\F_p^\times$, so $G_{K(\mu_p)}$ acts through $\{(a,a^{-1})\}$. Moreover $F'K$ is disjoint from $K(E[p])$, since it is ramified at primes dividing $D_{F'}$, which do not divide $pN_E$. Hence $\bar\chi_\fp/\bar\chi_\fpb=\bar\chi_\fp^2$ takes every value $a^2$, $a\in\F_p^\times$, on $G_{K(\mu_p)F'}$, and it is non-trivial because $p\ge5$.
\item Condition (vi), that $F'\ne\Q(\sqrt5)$ when $p=5$, is automatic, since $p$ is unramified in $F'$.
\end{itemize}

(2) The units of $\Lambda^{a,\rm ur}\otimes\Qp$ are $p^{\Z}\cdot(\Lambda^{a,\rm ur})^\times$, by the Weierstrass preparation theorem, so by (1) a generator $\mathcal F$ of $\mathrm{char}(X_{\rm Gr})$ satisfies $\mathcal F=p^ku\,L^\BDP_p(f)$ with $k\in\Z$ and $u$ a unit. By \cite[Thm.~5.1]{CHKLL}, whose hypotheses $E(F)[p]=0$ and the non-triviality of the localisations at $v$ and $\bar v$ hold here (Lemma~\ref{lem:consequences}(1) and $\lambda_a\ne0$), $\mathcal F$ lies in $I^2$ and its image in $I^2/I^3$ is, up to a $p$-adic unit, the right-hand side of \cite[Conj.~4.2]{AC21}; the order is exactly two because the cross pairing $\det h(x_i,\tau x_j)=\lambda_1^2\log_E(Q_v)^2$ is non-zero. By Corollary~\ref{cor:ACunit}, $\bar L^\BDP_p$ has the same valuation. Comparing leading coefficients gives $k=0$.
\end{proof}

\section{Numerical evidence}\label{sec:evidence}

This section collects the computations, in the order of the results they test: the conjecture itself through Theorems~\ref{thm:three} and~\ref{thm:digits} (\S\S\ref{ss:method}--\ref{ss:campaign}), the normalisations in rank one (\S\ref{ss:rankone}), the position of Rubin's class (\S\ref{ss:fourth-test}), the Heegner identity (\S\ref{ss:heegner-num}), the Birch and Swinnerton-Dyer quotient over $F$ (\S\ref{ss:bsdF-num}) and the two readings of the conjecture of Agboola and Castella (\S\ref{ss:readings}).

\subsection{\texorpdfstring{Computing $q_p$ and $q_\infty$}{Computing q\_p and q\_infty}}\label{ss:method}
By Corollary~\ref{prop:Cp} and Theorem~\ref{thm:digits}, $q_p$ modulo $p^{B_p(m)-2}$, which is $p^m$ for $p\ge m+3$, is determined by $m$ central values $M_k$ and by $j$. The central value $L(\psib^{\,2ke+1},ke+1)$ is computed by the approximate functional equation with $O(\sqrt{N_E}\,k p)$ terms and no cancellation, and recognised in $K$ (its denominator is a power of $p$); $j$ and $H$ are computed from \eqref{eq:hp} with $p$-adic division polynomials, in milliseconds, and $H$ was checked against \texttt{ellpadicregulator} of PARI/GP. The total cost is of the order of $\sqrt{N_E}\,p^3$ bit operations up to logarithmic factors, against linear algebra in a space of dimension proportional to $N_E$ for modular symbols. The algorithm, the certified recomputation in ball arithmetic of the values with $\Shan\ne1$, and the outputs of the computations are in the supplementary material.

On the complex side, $q_\infty$ is computed with the $L$-function routines of PARI/GP, and $\Shan=q_\infty\,\#E(\Q)_\tors^2/\Tam$ is recognised as an integer: to $100$ decimal digits for $E_{56}$, to $60$ for $E_{17},E_{-33},E_{-34},E_{-39}$, and for the $1368$ curves with $\Shan\in\{4,9\}$ of \S\ref{ss:fourth-test} to the working precision of $38$ digits, the largest distance $|\Shan-\mathrm{round}(\Shan)|$ being $3\cdot10^{-37}$. For curves of small conductor, $c_2(p)$ is also computed from modular symbols (SageMath/eclib), independently of the Katz measure; for instance $q_5/q_\infty=1+O(5^8)$ for $E_{56}$ and $q_{13}/q_\infty=1+O(13^3)$, with $q_\infty$ the rational number recognised above, and in all $22$ pairs computed in this way the two sides agree to the full precision.

\subsection{\texorpdfstring{The conjecture modulo $p^3$}{The conjecture modulo p\^{}3}}\label{ss:campaign}
We tested the curves $E_D$ of rank two, proved by $2$-descent with \texttt{ellrank}, in three sets: all curves with $|D|\le6000$; those with $6001\le|D|\le40000$ and $\Shan\ne1$; and those with $40001\le|D|\le249996$ and $\Shan\notin\{1,4\}$; at the non-anomalous primes $p\in\{5,13,17,29,37,41\}$ within a cost bound. Among the $61\,261$ rank-two curves with $6001\le|D|\le249996$, $\Shan$ takes only the values $1$, $4$ ($1808$ curves) and $9$ ($20$ curves).

\begin{table}[ht]
\caption{The $5941$ pairs $(E_D,p)$, $1898$ curves, conductors up to $3.9\cdot10^{12}$. A pair is counted under ``$\equiv q_\infty\,(p)$'' if Corollary~\ref{prop:Cp} gives $q_p\equiv q_\infty\pmod p$, and under ``$\equiv q_\infty\,(p^3)$'' if Theorem~\ref{thm:digits} with three central values ($p\ge13$) or four ($p=5$) gives $q_p\equiv q_\infty\pmod{p^3}$. For $j'\equiv0$ Corollary~\ref{prop:Cp} predicts $M'\equiv0$.}\label{tab:summary}
{\footnotesize\setlength{\tabcolsep}{2.5pt}\begin{tabular}{rrrrrrrr}
\toprule
$\Sha_{\mathrm{an}}$ & $p$ & pairs & $j'\not\equiv0$ & $\equiv q_\infty\,(p)$ & $\equiv q_\infty\,(p^3)$ & $j'\equiv0$ & $M'\equiv0$\\
\midrule
1 & 5 & 759 & 648 & 648 & 648 & 111 & 111\\
1 & 13 & 1495 & 1378 & 1378 & 1378 & 117 & 117\\
1 & 17 & 1530 & 1447 & 1447 & 1447 & 83 & 83\\
1 & 29 & 886 & 853 & 853 & 853 & 33 & 33\\
1 & 37 & 458 & 443 & 443 & 443 & 15 & 15\\
1 & 41 & 326 & 314 & 314 & 314 & 12 & 12\\
\midrule
4 & 5 & 140 & 116 & 116 & 116 & 24 & 24\\
4 & 13 & 208 & 192 & 192 & 192 & 16 & 16\\
4 & 17 & 75 & 71 & 71 & 71 & 4 & 4\\
4 & 29 & 5 & 5 & 5 & 5 & 0 & 0\\
\midrule
9 & 5 & 19 & 17 & 17 & 17 & 2 & 2\\
9 & 13 & 20 & 19 & 19 & 19 & 1 & 1\\
9 & 17 & 20 & 17 & 17 & 17 & 3 & 3\\
\midrule
\multicolumn{2}{r}{total} & 5941 & 5520 & 5520 & 5520 & 421 & 421\\
\bottomrule
\end{tabular}
}
\end{table}

Table~\ref{tab:summary} summarises the result: in all $5520$ pairs with $j'\not\equiv0$ one has $q_p\equiv q_\infty\pmod{p^3}$, and in all $421$ pairs with $j'\equiv0$ one has $M'\equiv0$; there is no exception and no numerical failure. For a sample we went further. For $E_{56}$, twelve central values give $q_5\equiv1\pmod{5^{12}}$, the precision promised by Theorem~\ref{thm:digits} for $m=12$, independently of the lift $y_0$ (the dependence on $y_0$ is $O(5^{12})$); the modular-symbol value $q_5=1+O(5^8)$, independent of the Katz measure, is limited by the precision of the modular-symbol computation. For the curves $E_{-2379}$ and $E_{-1513}$ with $\Shan=4$ we obtained $q_p\equiv q_\infty$ modulo $5^4$ and $13^3$.

The curves with $\Shan=4$ and $9$ are the point of the experiment. The integers $4$ and $9$ are prime to $p$, so they are invisible to any statement about characteristic ideals; only the exact formula predicts that $q_{13}\equiv9$ rather than $1$ or $4$, and, modulo $13^3$, rather than any of the $2027$ other units. These factors are arithmetic: for $\Shan=4$ the $2$-descent gives $\Sha(E/\Q)[2]\supseteq(\Z/2)^2$, and for the $20$ curves with $\Shan=9$ an explicit $3$-descent in the style of Schaefer and Stoll \cite{SchaeferStoll} gives $\Sha(E/\Q)[3]\supseteq(\Z/3)^2$ unconditionally.

\subsection{Calibration in rank one}\label{ss:rankone}
The same extraction applies in rank one, where the answer is known. If $E_D$ has rank one with generator $y$, then $\mathcal F$ vanishes to order one, $\mathcal F_1(a,b)=\tfrac12(a-b)\mathcal F_1(1,-1)$, and \eqref{eq:rubin-rank1} together with Lemma~\ref{lem:BK} gives, for $c_k$ solving $\sum_kc_kk^n=\delta_{n,1}$ $(1\le n\le m)$ and $y_0\equiv a_p^{-1}$,
\[\begin{split}\sum_{k=1}^mc_k\,y_0^{\,k}\,\iop\bigl(\mathcal N(ke,0)\bigr)&\equiv\tfrac12\,\iop(f_0\omega_E)\,e\,(1-\alpha^{-1})^2\,h^+(y)\,\frac{L'(E,1)}{\Omega_E\,\hat h(y)}\\&\qquad\pmod{p^{\min_{n>m}(n-v_p(n!))}}.\end{split}\]
We verified this for $(D,p)=(14,5)$, $(6,5)$, $(-9,5)$, $(-19,5)$, $(14,13)$, $(21,13)$, $(2,13)$, in each case at least to the predicted precision (for instance, for $(14,5)$ the quotient of the two sides is $\equiv1\pmod{5^7}$); in $(2,13)$, where $v_{13}(h^+(y))=4$, seven edge values are needed. Since all constants are those of \S\ref{sec:katz}, this checks the normalisations in a case where the exact formula is a theorem.

\subsection{The position of Rubin's class}\label{ss:fourth-test}
Both sides of \eqref{eq:fourth} are computable: $\kappa_1$ from the values $y_0^{a+b}\iop\mathcal N(ae-1,be-1)$ at the $m(m+3)/2$ points $(a,b)$ of smallest weight, by a two-variable Vandermonde system as in Theorem~\ref{thm:digits} ($\kappa_1$ is then known modulo $p^{\min_{n>m}(n-v_p(n!))}$, independently of the $p$-adic period), and $h^+(Q_\lambda)$ from $H$ and the formal logarithms. Put
\[\eta_p=\frac{\kappa_1/(p-1)}{\tfrac12\,\iop(\bar\omega_E/f_0)\,(1-\alpha/p)^2\,q_\infty\,h^+(Q_\lambda)},\]
so that Conjecture~\ref{conj:main} predicts $\eta_p=1$.

\begin{table}[ht]
\caption{The identity of Theorem~\ref{thm:fourth} with $q_\infty$: $\eta_p$. $m$ is the degree of the Vandermonde system; $^\dagger$anomalous prime ($a_5\equiv1\pmod5$).}\label{tab:fourth}
{\footnotesize\begin{tabular}{rrrrrrr}
\toprule
$D$ & $p$ & $\Sha_{\mathrm{an}}$ & $q_\infty$ & $m$ &$\kappa_1$ known mod & $v_p(\eta_p-1)$\\
\midrule
$56$ & $5$ & $1$ & $1$ & $6$ & $5^{6}$ & $5$\\
$-34$ & $5$ & $1$ & $1$ & $6$ & $5^{6}$ & $5$\\
$-39$ & $5$ & $1$ & $2$ & $6$ & $5^{6}$ & $5$\\
$-25721$ & $5$ & $4$ & $8$ & $2$ & $5^{3}$ & $2$\\
$17$ & $13$ & $1$ & $1$ & $3$ & $13^{4}$ & $3$\\
$56$ & $13$ & $1$ & $1$ & $3$ & $13^{4}$ & $3$\\
$-33$ & $13$ & $1$ & $1$ & $3$ & $13^{4}$ & $3$\\
$-34$ & $13$ & $1$ & $1$ & $3$ & $13^{4}$ & $3$\\
$56$ & $17$ & $1$ & $1$ & $2$ & $17^{3}$ & $2$\\
\midrule
$17$ & $5^{\dagger}$ & $1$ & $1$ & $6$ & $5^{6}$ & $5$\\
$-33$ & $5^{\dagger}$ & $1$ & $1$ & $6$ & $5^{6}$ & $5$\\
\bottomrule
\end{tabular}
}
\end{table}

In every row of Table~\ref{tab:fourth} one has $\eta_p=1$ to the full available precision: five $5$-adic digits for $E_{56}$, $E_{-34}$ and $E_{-39}$, three $13$-adic digits for four curves, and two $5$-adic digits for $E_{-25721}$, with $\Shan=4$ and $q_\infty=8$, where $\Shan=1$, i.e.\ $q_\infty=2$, would give $\eta_5\equiv4\pmod{25}$. Any other fourth root of unity in \eqref{eq:tau}, or any stray power of $p$ in Lemma~\ref{lem:rubin-norm}, would be visible in the first digit. The identity also holds at the two anomalous pairs $(17,5)$ and $(-33,5)$, where $\det H$ is a unit.
We then ran the test at $p=5$ on the curves with $\FourthWhich$, in order of conductor and with $m\in\{\FourthKs\}$: $\FourthDone$ curves ($\FourthShaFour$ with $\Shan=4$, $\FourthShaNine$ with $\Shan=9$, $\FourthAnom$ of them anomalous) with conductors from $\FourthNmin$ to $\FourthNmax$ and $q_\infty\in\FourthAE$. Of these, $\FourthVthree$ give $\eta_5\equiv1\pmod{5^3}$, $\FourthVtwo$ modulo $5^2$, $\FourthVone$ modulo $5$, and $\FourthVzero$ have no digit available because of a large valuation of $\kappa_1$; there is no discrepancy in any digit computed.

\subsection{Derived Heegner classes}\label{ss:heegner-num}
We verified \eqref{eq:heegner0} for two pairs, computing $\zeta^0$ modulo a power of $p$ from Heegner points.

\emph{Method.} For conductor $p^{s+1}$, the Heegner points on $X_0(N_E)$ are computed as elliptic logarithms $z=\sum a_n n^{-1}q^n$ at CM points $\tau$, one per class of $\mathrm{Pic}(\Z+p^{s+1}\OK_F)$, choosing among the four Atkin--Lehner families the representative of largest imaginary part (the families are related by explicit affine maps on the class group). The sums over the cosets of the subgroup of index $p^s$ are the conjugates of a point of $E(F_s)$, where $F_s$ is the layer of degree $p^s$ of $F^a_\infty$ over $F$; the $x$-coordinates of these conjugates are recognised as the roots of a polynomial in $\Q[x]$ of degree $p^s$, the Galois group of $F_s/\Q$ being dihedral. Then $\tilde y_s=\alpha^{-s-1}S_{s+1}-\alpha^{-s-2}S_s$, with $S_s$ the coset sums of conductor $p^s$, and $\mathcal D\equiv-\sum_i i\,\sigma^i\tilde y_s\pmod{p^sE(F)}$. This is evaluated in $E(\F_\ell)$ for primes $\ell$ that split completely in $F_s$ and for which $E(\F_\ell)$ has a cyclic factor of order $p^s$: the cyclic order of $\sigma$ on the roots modulo $\ell$ is fixed by rational invariants of the polynomial, the signs of the $y$-coordinates by the trace relation and the Euler system relation $\mathrm{Norm}_{F_{s}/F_{s-1}}(S_{s+1})=a_pS_s$ up to torsion, and $\mathcal D$ is then compared with $Q_v$ in $E(\F_\ell)/p^s$. Finally $\ell_a(\sigma_a)$ is computed as $(\log_p\pi-\log_p\bar\pi)/h_F$ for a generator $\pi$ of $\mathfrak q^{h_F}$, $\mathfrak q$ a prime whose class has image $\sigma_a^{\pm1}$.

\emph{$E_{56}$, $F=\Q(\sqrt{-31})$, $p=5$.} Here $s_F=4$, $h_F=3$, $(1-\alpha^{-1})^4\equiv1\pmod{25}$. At conductor $25$ ($60$ points, $1000$ digits) the coset sums are the roots of an irreducible quintic in $\Q[x]$ whose denominators are squares of primes, one of them with $200$ digits; at conductor $125$ ($300$ points, $8000$ digits) of an irreducible polynomial of degree $25$ with coefficients of about $4900$ bits, recognised far below the precision limit, the Euler system relation holding to $10^{-1001}$. Eight split primes $\ell$ give $\zeta^0$ modulo $25$ consistently, with $\mathcal D$ on the predicted line $Q_v$, and the right-hand side of \eqref{eq:heegner0} is $\equiv1\pmod{25}$, while $q_5=1+O(5^{12})$ (\S\ref{ss:campaign}). Since $q_5$ is known independently, this confirms \eqref{eq:heegner0}, sign included.

\begin{table}[ht]
\caption{The Heegner identity for $E_{56}$ over $F=\Q(\sqrt{-31})$ at $p=5$. Here $\Xi=2\ell_a(\sigma_a)\zeta^0/(1-\alpha^{-1})^2$, so that \eqref{eq:heegner0} reads $q_p=\Xi^2/s_F$.}\label{tab:heegner}
{\footnotesize\begin{tabular}{ll}
\toprule
$h_F$; $\mathrm{Pic}$ of conductor $25$, $125$ & $3$; $\Z/60$, $\Z/300$\\
$q_\infty(E^\chi)$, $c_E^{-2}\kappa_F\deg\varphi$, $s_F$ & $4$, $2$, $4$\\
$a_5$, $\alpha\bmod25$, $(1-\alpha^{-1})^4\bmod125$ & $-2$, $13$, $76$\\
$\mathcal D\bmod5$ (conductor $25$) & $\pm(2P_1+3P_2)$, on the line of $Q_v$\\
$\Xi\bmod25$ (conductor $125$) & $\pm2$\\
$\Xi^2/s_F\bmod25$ & $1$\\
$q_5$ (Theorem~\ref{thm:digits}; modular symbols) & $1+O(5^{12})$; $1+O(5^8)$\\
\bottomrule
\end{tabular}}
\end{table}

\emph{$E_{-289}$, $F=\Q(\sqrt{-111})$, $p=5$.} Here $h_F=8$, $\mathrm{Pic}$ of conductor $25$ is $\Z/80\times\Z/2$, $s_F=4$ (with $q_\infty(E^\chi)=8$ and $c_E^{-2}\kappa_F\deg\varphi=1$, against $4$ and $2$ for $E_{56}$). This case was computed after the constant of \eqref{eq:heegner0} had been fixed. At conductor $25$ ($160$ points, $4000$ digits) the quintic has coefficients of about $3400$ bits; thirty split primes give $\mathcal D\equiv\pm2P_2\pmod5$ in the basis $P_1,P_2,P_3$, on the predicted line ($Q_v\equiv(0,1,0)$), and \eqref{eq:heegner0} gives $q_5\equiv1\pmod5$, which is the known value.

The data and scripts are in the supplementary material.

\subsection{\texorpdfstring{The Birch and Swinnerton-Dyer quotient over $F$}{The Birch and Swinnerton-Dyer quotient over F}}\label{ss:bsdF-num}
\begin{table}[ht]
\caption{Terms of \eqref{eq:bsdF} for a sample of the $79$ pairs $(E_D,\Q(\sqrt{-d}))$; $\#E(F)_\tors=2$ in all cases. In every row $\Tam s_F/\#E(\Q)^2_\tors$ is a square.}\label{tab:bsdF}
{\footnotesize\begin{tabular}{rrccccccc}
\toprule
$D$ & $d$ & $q_\infty(E)$ & $q_\infty(E^\chi)$ & $c_E^{-2}\kappa_F\deg\varphi$ & $s_F$ & $j$ & $\prod_vc_v(E/F)$ & $\Shan(E/F)$\\
\midrule
$56$ & $31$ & $1$ & $4$ & $2$ & $4$ & $1$ & $16$ & $1$\\
$56$ & $55$ & $1$ & $16$ & $2$ & $16$ & $0$ & $16$ & $1$\\
$56$ & $223$ & $1$ & $16$ & $2$ & $16$ & $1$ & $16$ & $4$\\
$17$ & $111$ & $1$ & $4$ & $2$ & $4$ & $1$ & $16$ & $1$\\
$-68$ & $111$ & $2$ & $16$ & $1$ & $8$ & $1$ & $64$ & $1$\\
$-289$ & $111$ & $1$ & $8$ & $1$ & $4$ & $1$ & $16$ & $1$\\
$-999$ & $71$ & $2$ & $16$ & $1$ & $8$ & $1$ & $64$ & $1$\\
$-594$ & $95$ & $2$ & $16$ & $1$ & $8$ & $1$ & $64$ & $1$\\
$-1224$ & $191$ & $2$ & $16$ & $1$ & $8$ & $1$ & $64$ & $1$\\
$-89$ & $39$ & $1/2$ & $4$ & $1$ & $2$ & $1$ & $4$ & $1$\\
$-73$ & $183$ & $1/2$ & $16$ & $1$ & $8$ & $1$ & $4$ & $4$\\
$-113$ & $31$ & $1/2$ & $4$ & $1$ & $2$ & $1$ & $4$ & $1$\\
\bottomrule
\end{tabular}
}
\end{table}

We computed the terms of \eqref{eq:bsdF} for $79$ pairs $(E,F)$ (Table~\ref{tab:bsdF}), with ten curves
\[D\in\{56,17,-68,-289,-999,-594,-89,-73,-113,-1224\}\]
and all admissible $F=\Q(\sqrt{-d})$ with $d\le260$ for which a generator of $E^\chi(\Q)$ is found; the index $2^j$ is determined by testing $2$-divisibility in $E(F)$ of the combinations of $P_1,P_2,P_3$ and torsion. In all cases $j\in\{0,1\}$, $\#E(F)_\tors=2$, $\prod_vc_v(E/F)\in\{4,16,64\}$, $c_E^{-2}\kappa_F\deg\varphi=2$ for $D>0$ and $1$ for $D<0$ (Lemma~\ref{lem:kappa}), and $\Shan(E/F)$ is an integer to $40$ digits, equal to $1$ in $75$ cases and to $4$ in $4$; in particular $\Tam\,s_F/\#E(\Q)^2_\tors\in\{1,4,16\}$ is a square in every case, as Theorem~\ref{thm:tamsquare} predicts; indeed $\Tam(E)s_F\in\{4,16,64\}$. The test is not vacuous: for the seven curves with $q_\infty(E)\in\{2,\tfrac12\}$, such as $E_{-68}$, which is $2$-isogenous to $E_{17}$, a square $s_F$ would contradict Conjecture~\ref{conj:main} at $p=5$, where $2$ is not a square; one finds instead $c_E^{-2}\kappa_F\deg\varphi=1$ and $q_\infty(E^\chi)\in\{4,16\}$, so that $s_F\in\{2,8\}$ and $q_\infty s_F$ is a square.

\subsection{The two readings of the conjecture of Agboola and Castella}\label{ss:readings}
By \S\ref{sec:AC}, Conjecture~\ref{conj:AC} with the sesquilinear product is equivalent to $\Xi^2=q_\infty s_F$, and with the ring product to $\Xi^2=-q_\infty s_F$; these statements involve only Heegner points and the complex side. For $(E_{56},\Q(\sqrt{-31}),5)$ the Heegner points give $\Xi^2\equiv4\pmod{25}$ (Table~\ref{tab:heegner}), and for $(E_{-289},\Q(\sqrt{-111}),5)$ they give $\Xi^2\equiv4\pmod5$; in both cases $q_\infty s_F=4$, whereas $-4\equiv21\pmod{25}$ and $-4\equiv1\pmod5$. Both examples therefore agree with the sesquilinear reading and not with the ring product.

\section{Concluding remarks}\label{sec:final}

\subsection{What is proved and what is not}
Theorems~\ref{thm:three}, \ref{thm:digits}, \ref{thm:position}, \ref{thm:fourth} and \ref{thm:heegner} are identities between $p$-adic quantities in which $q_p$ appears as the only unknown: $q_p$ is the constant of the quadratic term of the Katz measure, and the coordinate of Rubin's elliptic-unit class and of the derived Heegner class. Together with \eqref{eq:rubin-rank1} for the twist they reduce Conjecture~\ref{conj:main} to any one of several explicit statements: the position of Rubin's class, the value of one derivative of the Katz measure at the conjugate character, the coordinate $\zeta$ of the derived Heegner class, or the leading coefficient of the BDP $p$-adic $L$-function. Each of them is an identity between a $p$-adic number and the complex number $q_\infty$. What can be proved without knowing $q_\infty$ is what does not see $q_p$ beyond its valuation and its square class (for the unit $q_p/A_E$, when $\Sha(E/\Q)$ is finite): the square class itself (Theorem~\ref{thm:squares}, whenever an auxiliary field $F$ as in \S\ref{ss:hyp} exists) and, combined with the known $p$-parts of the Birch and Swinnerton-Dyer formula, the conjecture of Agboola and Castella up to a $p$-adic unit (Corollary~\ref{cor:ACunit}), and with \cite{BCS25} the integral anticyclotomic main conjecture (Theorem~\ref{thm:bdpMC}). In rank one, the corresponding identity is proved by computing one global object, the Heegner point, both complex-analytically and $p$-adically. In rank two the $p$-adic objects that we control (Rubin's class, the derived Heegner class, Kato's derived class) lie on lines of $E(\Q)\otimes\Qp$ that are not defined over $\Q$, and no archimedean formula for them is known. Every identity available to us is homogeneous in $q_p$, and its square class comes from the quadratic nature of heights and of the Cassels--Tate pairing.

\subsection{Limitations}
(1) The values of the Katz measure are recognised from numerical approximations of Hecke $L$-values; the decisive ones were recomputed in certified ball arithmetic, conditional on an a priori bound for their denominators (supplementary material).
(2) Only $K=\Q(i)$ is treated; for the other class number one fields except $\Q(\sqrt{-3})$, where $s_2=0$ as well, $s_2$ need not vanish and the reduced theta function replaces $\sigma$.
(3) The equivalence of Theorem~\ref{thm:AC}(1) with Conjecture~\ref{conj:main} reads the products of anticyclotomic heights sesquilinearly, a choice supported numerically by both examples (\S\ref{ss:readings}) and by the comparison with \cite[Conj.~4.3]{AC21}.
(4) The identity \eqref{eq:LfLK} enters only through \eqref{eq:ratio}, which does not see its constant. The normalisation \eqref{eq:T3} is taken from Perrin-Riou's factorisation \cite[(1.1)]{PR87}, and the constants of the explicit reciprocity law from \cite[Thm.~5.7]{CastellaHsieh18} and \cite[Thm.~3.7]{Castella20}. The numerical verification of \S\ref{ss:heegner-num} checks the resulting identity modulo $25$ for one pair and modulo $5$ for another. Theorem~\ref{thm:AC} and Corollary~\ref{cor:ACunit} depend on the same constants.

\subsection{Questions}
\begin{enumerate}
\item Is there an archimedean formula for the derived Heegner class $\mathcal D$, or for Rubin's class, in the spirit of Gross--Zagier? By Theorem~\ref{thm:fourth}, when $h^+(Q_\lambda)\ne0$, it would suffice to prove \eqref{eq:fourth} with $q_\infty$: one derivative of the Katz measure at $\psi^*$ against one cyclotomic height.
\item The coordinate $\zeta^0$ of $\mathcal D^0$ is, up to the explicit factor $\ell_a(\sigma_a)/(1-\alpha^{-1})^2$, predicted by Conjecture~\ref{conj:AC} to be $\pm2^{-j}\sqrt{\Shan(E/F)}\prod_\ell c_\ell(E)/\#E(F)_\tors$; this is rational if moreover $\Shan(E/F)$ is a rational square, as the Birch and Swinnerton-Dyer conjecture and the Cassels--Tate pairing predict, which is more than $q_p=q_\infty$. Can its rationality be proved directly, for instance by $p$-adic deformation along the Hida family of CM forms $\theta(\psi^{k-1})$, in which, in all the weights we computed ($k\le22$ for $E_{56}$), $k=2$ is the only one where $L(f_k/F,s)$ vanishes to order three, the others having order one?
\end{enumerate}

\bibliographystyle{amsplain}
\bibliography{refs}

\end{document}